\documentclass[10pt,a4paper]{article}

\usepackage{graphicx} 
\usepackage{amsmath}
\usepackage{amsfonts}
\usepackage{threeparttable,url}
\usepackage{multirow}
\usepackage{algorithm}
\usepackage{algorithmic}
\usepackage{color}
\usepackage[font=footnotesize,labelfont=bf]{caption}

\usepackage{cite}

\newtheorem{remark}{Remark}

\usepackage[top=1.25in,right=1.25in,left=1.25in,bottom=1.25in]{geometry}

\begin{document}

%
%


\title{Large-scale optimization with the primal-dual column generation method}
\author{Jacek Gondzio\thanks{School of Mathematics, University of Edinburgh, United Kingdom (j.gondzio@ed.ac.uk.)} \ \ \ \ 
Pablo Gonz\'{a}lez-Brevis\thanks{School of Engineering, Universidad del Desarrollo, Concepci{\'o}n, Chile (pablogonzalez@ingenieros.udd.cl)}
\ \ \ \ Pedro Munari\thanks{Production Engineering Department, Federal University of S\~ao Carlos, Brazil (munari@dep.ufscar.br)}}

\date{}
\maketitle
\vspace{-3ex}
\begin{center}
School of Mathematics, University of Edinburgh\\
The King's Buildings, Edinburgh, EH9 3JZ, UK\\
Technical Report ERGO 13-014, \\
\date{January 22, 2015}.
\end{center}

\maketitle

\allowdisplaybreaks

\begin{abstract}
The primal-dual column generation method (PDCGM) is a general-purpose column generation technique that relies on the primal-dual interior point method to solve the restricted master problems. The use of this interior point method variant allows to obtain suboptimal and well-centered dual solutions which naturally stabilizes the column generation. 
{As recently presented in the literature,}
reductions in the number of calls to the oracle and in the CPU times are typically observed when compared to the standard column generation, which relies on extreme optimal dual solutions. 
{However, these results are based on relatively small problems obtained from linear relaxations of combinatorial applications.}
In this paper, we investigate the behaviour of the PDCGM {in a broader context, namely} when solving large-scale convex optimization problems. We have selected applications that arise in important real-life contexts such as data analysis (multiple kernel learning problem), decision-making under uncertainty (two-stage stochastic programming problems) and telecommunication and transportation networks (multicommodity network flow problem). In the numerical experiments, we use publicly available benchmark instances to compare the performance of the PDCGM against
{recent results for different methods presented in the literature, which were the best available results to date.
The analysis of these results suggests} that the PDCGM offers an attractive alternative over specialized methods since it remains competitive in terms of number of iterations and CPU times { even for large-scale optimization problems.}

\vspace{3ex}

\end{abstract}

\section{Introduction}\label{intro}

Column generation is an iterative oracle-based approach which has been widely used in the context of continuous as well as discrete optimization \cite{BriLemMeuMicPerVan08,LubDes05}. In this method, an optimization problem with a huge number of variables is solved by means of a reduced version of it, the restricted master problem (RMP). At each iteration, the RMP is modified by the addition of columns which are generated by the oracle (or, pricing subproblem). To generate these columns, the oracle uses {a} dual solution of the RMP. 

In the standard column generation, {optimal dual solutions of the RMP are used in the oracle to generate new columns. Since a simplex method is typically used to optimize the RMP, these solutions correspond to extreme points of the dual feasible set of the RMP. As a result,} large variations are typically observed between dual solutions of consecutive {column generation} iterations, a behavior that may cause the slow convergence of the method. In addition, when active-set methods, such as the simplex method, are used to solve the RMP, degeneracy {may} adversely affect the performance of the column generation method. These drawbacks {are also observed} in the cutting plane method \cite{Kel60}, which is the dual counterpart of column generation. Several alternatives to overcome such weaknesses have been proposed in the literature. Some of them modify the RMP by adding penalty terms and/or constraints to it with the purpose of limiting the large variation of the dual solutions \cite{Marsten1975,duMerle1999,BriLemMeuMicPerVan08,Frangioni2002}. {Other alternatives use dual price smoothing techniques \cite{wentges1997,neame2000nonsmooth}}. Finally, there exist variants of column generation which rely on naturally stable approaches to solve the RMP, such as interior point methods \cite{goffin1992,mitchell1996,martinson1999,GonGonMun2013}.

The primal-dual column generation method (PDCGM) \cite{GonGonMun2013} is a variant which relies on well-centered and suboptimal dual solutions of the RMPs. To obtain such solutions, the method uses the primal-dual interior point algorithm \cite{Gondzio2011}. The optimality tolerance used to solve the restricted master problems is loose during the first iterations and it is dynamically reduced as the method approaches optimality.
This reduction guarantees that an optimal solution of the original problem is obtained at the end of the procedure. Encouraging computational results are reported in \cite{GonGonMun2013} regarding the use of the PDCGM to solve the relaxations of three widely studied mixed-integer programming problems, namely the cutting stock problem, the vehicle routing problem with time windows and the capacitated lot-sizing problem with setup times, after applying a Dantzig-Wolfe decomposition \cite{DanWol61} to their standard compact formulations.

In this paper, we extend the computational study presented in \cite{GonGonMun2013} by analyzing the performance of the PDCGM applied to solve large-scale convex optimization problems.
{The applications considered in \cite{GonGonMun2013} have relatively small restricted master problems and the bottleneck is 
{ in solving the oracle.}
On the other hand, the applications we address in the current paper have large restricted master problems and the oracle subproblems are relatively easy to solve. Hence, we evaluate the performance of PDCGM operating in very different conditions besides addressing a broader class of optimization problems.}

By large-scale problems we mean a formulation which challenges the current state-of-the-art implementations of optimization methods, due to a very large number of constraints and/or variables. Furthermore, we assume {that such} formulation has a special structure which allows the use of a reformulation technique, such as the Dantzig-Wolfe decomposition. Hence, large-scale refers not only to size, but also structure.
The problems we address arise in important real-life contexts such as data analysis, decision-making under uncertainty and telecommunication/transportation networks. 

The main contributions of this paper are the following. First, we review three applications which have gained a lot of attention in the optimization community in the past years and describe them in the column generation framework. Second, we study the behavior of the PDCGM to solve publicly available instances of these applications and compare its performance with 
{recent results of other}
stabilized column generation/cutting plane methods,
{which we believe are the best results presented in the literature for the addressed problems.}
{As a third contribution}, we make our software available for any research use.

The remainder of this paper is organized as follows. 
In Section \ref{sec:dwdcolgen}, we describe the decomposition principle and the column generation technique for different situations.
In Section \ref{sec:algo}, we outline some stabilized column generation and  cutting plane algorithms which have proven to be effective for the problems we deal with in this paper. 
In Sections \ref{sec:mkl}, \ref{sec:stoch} and \ref{sec:nlmcnf} we describe the \emph{multiple kernel learning} (MKL), the \emph{two-stage stochastic programming} (TSSP) and the \emph{multicommodity network flow} (MCNF) problems, respectively. 
In 
{each
of these three} sections, we present the problem formulation and derive the column generation components, namely the master problem and oracle.
Then, we report on computational experiments comparing the performance of PDCGM with other state-of-the-art techniques. 
Finally, we summarize the main outcomes of this paper in Section \ref{sec:conclusions}.
\section{Reformulations and the column generation method}\label{sec:dwdcolgen}

Consider an optimization problem stated in the following form
\begin{equation}
\label{eq:original}
\min \ \ c^T x, \ \mbox{s.t.} \ \ Ax \leq b, \ \ \ x \in \mathcal{X} \subseteq \mathbb{R}^n,
\end{equation}
where {$x$ is a vector of decision variables, $\mathcal{X}$ is a non-empty convex set and $c \in \mathbb{R}^n$, $A \in \mathbb{R}^{m \times n}$ and $b \in \mathbb{R}^m$ are the problem parameters}. We assume that without the presence of the set of \textit{linking} constraints $Ax \leq b$, problem \eqref{eq:original} would be easily solved by taking advantage of the structure of $\mathcal{X}$. {More specifically, the $n$ variable indices can be suitably partitioned into $K$ subsets, such that $\mathcal{X}$ is given by the Cartesian product $\mathcal{X}:= \mathcal{X}^1\times \dots \times \mathcal{X}^{K}$, in which the sets $\mathcal{X}^{k}$ are independent from each other. Following this notation, we have the partition $x = (x^{1}, \ldots, x^{K})$ with $x^k \in \mathcal{X}^k$ for every $k \in  \mathcal{K} = \{1, \ldots, K\}$.}

{We further assume that each set $\mathcal{X}^k$ can be described as a polyhedral set, either equivalently or by using a fine approximation as discussed 
{ later}
in Section \ref{sub:sec:convex:function}.} Hence, let ${P}^k$ and ${R}^k$ {denote} the sets of indices of all the extreme points and extreme rays of $\mathcal{X}^k$, respectively. With this notation, we represent by $x^{k}_p$ and $x^{k}_r$ the extreme points and the extreme rays of $\mathcal{X}^k$, with $p \in {P}^k$ and $r \in {R}^k$.
Any point $x^k \in \mathcal{X}^k$ can be represented as a combination of these extreme points and extreme rays as follows
\begin{equation}
\label{eq:rewritex}
x^k = \sum_{p \in {P}^k}  \lambda^k_p x^k_p + \sum_{r \in {R}^k} \lambda^k_r x^k_r, \mbox{ with } \sum_{p \in {P}^k}\lambda^k_p = 1,
\end{equation}
with the coefficients $\lambda^k_p \geq 0$ and $\lambda^k_r \geq 0$, for $p \in {P}^k$ and $r \in {R}^k$.
The Dantzig-Wolfe decomposition principle (DWD) \cite{DanWol61} consists in using this relationship to rewrite the original variable vector $x$ in problem (\ref{eq:original}).
Additionally, considering the partitions $c = (c^1, \ldots, c^{K})$ and $A = (A^1, \ldots, A^{K})$ which are induced by the structure of $\mathcal{X}$, we define , $c^k_p:= (c^k)^T x^k_p$ and $a^k_p:=A^k x^k_p$ for every $k \in \mathcal{K}$, $p \in {P}^k$, and $c^k_r:= (c^k)^T x^k_r$ and $a^k_r := A^k x^k_r$ for every $k \in \mathcal{K}$, $r \in {R}^k$. By using this notation, we can rewrite \eqref{eq:original} in the following equivalent form, known as the master problem (MP) formulation
\begin{eqnarray}
\min_{\lambda} & \displaystyle \sum_{k \in \mathcal{K}}\left(\sum_{p \in {P}^k}c_p^k\lambda_p^k + \sum_{r \in {R}^k}c_r^k\lambda_r^k\right), \label{eq:disagg:obj} \\
\mbox{s.t.} & \displaystyle \  \sum_{k \in \mathcal{K}}\left(\sum_{p \in {P}^k}a_p^k\lambda_p^k + \sum_{r \in R^k}a_r^k\lambda_r^k\right) \leq b, & \label{eq:disagg:c1}\\
& \displaystyle \ \sum_{p \in {P}^k}\lambda_p^k = 1, &\forall k \in \mathcal{K}, \label{eq:disagg:c2} \\
& \displaystyle \ \lambda_p^k \geq 0, &\forall k \in \mathcal{K},\forall p \in {P}^k, \label{eq:disagg:c3} \\
& \displaystyle \ \lambda_r^k \geq 0, &\forall k \in \mathcal{K},\forall r \in {R}^k. \label{eq:disagg:c4}
\end{eqnarray}
Notice that the coefficients in (\ref{eq:rewritex}) are now the decision variables in this model.
Constraints (\ref{eq:disagg:c1}) are called \textit{linking} constraints as they correspond to the set of constraints $Ax \leq b$ in (\ref{eq:original}). The \textit{convexity} constraints (\ref{eq:disagg:c2}) ensure the convex combination required by (\ref{eq:rewritex}), for each $k \in \mathcal{K}$. 
%
{As a result,} the value of the optimal solution of (\ref{eq:disagg:obj})-(\ref{eq:disagg:c4}), denoted by $z_{MP}$, is the same as the optimal value of (\ref{eq:original}).

\subsection{Column generation}

The number of extreme points and extreme rays in the MP formulation (\ref{eq:disagg:obj})-(\ref{eq:disagg:c4}) may be excessively large. Therefore, solving this problem by a direct approach is practically impossible. Moreover, in many occasions the extreme points and extreme rays which describe $\mathcal{X}$ are not available and have to be generated by a procedure which might be costly. Hence, we rely on the column generation technique \cite{ford1958,gilmore1961,LubDes05}, which is an iterative process that works as follows. At each \textit{outer} iteration, we solve a \textit{restricted master problem} (RMP), which has only a small subset of the columns/{variables} of the MP. The dual solution of the RMP is then used in the \textit{oracle} {with the aim of generating} one or more new extreme points or extreme rays, which {may} lead to new 
{columns}. 

Then, {if at least one of these new columns has a negative reduced cost, we add it to the RMP and start a new iteration. This step is called an outer iteration, to differentiate to the inner iterations which are the ones required to solve the RMP.}
{ The iterative process terminates when we can guarantee that the optimal solution of the current RMP is also optimal for the MP, even though not all the columns have been generated.}

The RMP has the same formulation as the MP, except that only a subset of extreme points and extreme rays are considered. In other words, the RMP is defined by the subsets $\overline{P}^k \subseteq P^k$ and $\overline{R}^k \subseteq R^k$, $k \in \mathcal{K}$. 
These subsets may change at every outer iteration by adding/removing extreme points and/or extreme rays. {Therefore,} the corresponding RMP can be represented as
\begin{eqnarray}
\min_{\lambda} & \displaystyle \sum_{k \in \mathcal{K}}\left(\sum_{p \in \overline{P}^k}c_p^k\lambda_p^k + \sum_{r \in \overline{R}^k}c_r^k\lambda_r^k\right), \label{eq:rmp:obj} \\
\mbox{s.t.} & \displaystyle \  \sum_{k \in \mathcal{K}}\left(\sum_{p \in \overline{P}^k}a_p^k\lambda_p^k + \sum_{r \in \overline{R}^k}a_r^k\lambda_r^k\right) \leq b, & \label{eq:rmp:c1}\\
& \displaystyle \ \sum_{p \in \overline{P}^k}\lambda_p^k = 1, &\forall k \in \mathcal{K}, \label{eq:rmp:c2} \\
& \displaystyle \ \lambda_p^k \geq 0, &\forall k \in \mathcal{K},\forall p \in \overline{P}^k, \label{eq:rmp:c3} \\
& \displaystyle \ \lambda_r^k \geq 0, &\forall k \in \mathcal{K},\forall r \in \overline{R}^k. \label{eq:rmp:c4}
\end{eqnarray}
The value of a feasible solution of the RMP (${z}^{ }_{RMP}$) provides an upper bound of the optimal value of the MP ($z_{MP}$), as this solution corresponds to a feasible solution of the MP in which all the components in ${P}^k \setminus \overline{P}^k$ and ${R}^k \setminus \overline{R}^k$ are equal to zero.

To verify if the current optimal solution of the RMP is also optimal for the MP, we use the dual solutions. Let $\overline{u} \in \mathbb{R}^{m}_{-}$ and $\overline{v} \in \mathbb{R}^{K}$ represent the dual optimal solutions associated with constraints \eqref{eq:rmp:c1} and \eqref{eq:rmp:c2}, respectively. The oracle checks the dual feasibility of these solutions in the MP by means of the reduced cost information. However, instead of calculating the reduced costs for every single {variable} of the MP, the oracle solves a pricing subproblem of the type
\begin{equation}
\label{sub:prob}
{SP}^k(\overline{u}) := \min_{x^k \in \mathcal{X}^k} \ \{ (c^k - \overline{u}^T A^k)^T x^k \},
\end{equation}
for each $k \in K$. There are two possible cases: 
either subproblem \eqref{sub:prob} has a bounded optimal solution, or its solution is unbounded.
In the first case, the optimal solution corresponds to an extreme point $x^k_p$ of $\mathcal{X}^k$. Let ${z}^{k}_{SP}(\overline{u}, \overline{v}) = (c^k - \overline{u}^T A^k)^T x^k_p - \overline{v}^k$ be the reduced cost of the variable associated to the column which is generated by using $x^k_p$. 
If this reduced cost is negative, then $p \in {P}^k \setminus \overline{P}^k$ and $x^k_p$ 
defines a new 
{column that should be added to the RMP.}
On the other hand, if ${SP}^k(\overline{u})$ has an unbounded solution, then we take the direction of unboundedness as an extreme ray $x^k_r$ for this problem.
The reduced cost in this case is given by ${z}^{k}_{SP}(\overline{u}, \overline{v}) = (c^k - \overline{u}^T A^k)^T x^k_r$.
If this reduced cost is negative, $x^k_r$ defines a
column which must be added to the RMP. By summing over all the negative reduced costs, we define the \textit{value of the oracle} as
\begin{equation}
{z}^{ }_{SP}(\overline{u}, \overline{v}) := \sum_{k \in \mathcal{K}} \min \{0; {z}^{k}_{SP} (\overline{u}, \overline{v}) \}.
\end{equation}
By using this value, we obtain the following relationship at any outer iteration of the column generation method
\begin{equation}
{z}^{ }_{RMP} + {z}^{ }_{SP}(\overline{u}, \overline{v}) \leq z_{MP} \leq {z}^{ }_{RMP}.
\end{equation}
When ${z}^{ }_{SP}(\overline{u}, \overline{v}) = 0$, dual feasibility has been achieved and hence the optimal solution of the RMP is also an optimal solution of the MP.

The standard column generation based on {optimal} dual solutions as we have just described {is known to have several drawbacks, specially when the simplex method is used to optimize the RMPs. In such case the resulting dual solutions correspond to extreme points of the dual feasible set. This typically leads to} large oscillations of consecutive dual solutions (\textit{bang-bang effect}), in particular at the beginning of the iterative process (\textit{heading-in effect}) \cite{vanderbeck2005, LubDes05}. Moreover, optimal dual solutions contribute also to the tailing-off effect, a slow progress of the method close to {termination. All these} drawbacks may result in a {large} number of calls to the oracle, slowing down the column generation method (similar issues can be observed in a cutting plane context). Therefore, several stabilization techniques have been proposed to overcome these limitations. In Section \ref{sec:algo}, we recall some of the techniques which have been successfully used within the applications {that} we address in this paper.

\subsection{Aggregated formulation}

{Formulation (\ref{eq:disagg:obj})-(\ref{eq:disagg:c4}) is often called the \emph{disaggregated} master problem, as each column in that formulation comes from an extreme point or an extreme ray of only one set $\mathcal{X}^k$, $k \in \mathcal{K}$. The separability of the sets $\mathcal{X}^k$ is reflected in the master problem formulation, as we have master variables $\lambda^{k}$ associated to each set $\mathcal{X}^{k}$, $k \in \mathcal{K}$, and $K$ convexity constraints in (\ref{eq:disagg:c2}). In some situations, it may be interesting to keep the variables aggregated in the master problem, e.g. when the number $K$ is too large so that the disaggregated formulation has a relatively large number of constraints and columns, as pointed out in \cite{ruszczynsky1986}. Hence, we can use the \textit{aggregated formulation} of the master problem, in which we do not exploit the structure of $\mathcal{X}$ explicitly when defining the master variables.} More specifically, let $P$ and $R$ be the sets of indices of extreme points and extreme rays of $\mathcal{X}$, respectively. Any given point $x \in \mathcal{X}$ can be written as the following combination
\begin{equation}
\label{eq:rewritexagg}
x = \sum_{p \in {P}} \lambda_p x_p + \sum_{r \in {R}} \lambda_r x_r, \mbox{ with } \sum_{p \in {P}} \lambda_p = 1,
\end{equation}
where $x_p$ and $x_r$ are extreme points and extreme rays of $\mathcal{X}$. 
By replacing this relationship in problem (\ref{eq:original}) and using the notation $c_p := c^T x_p$ and $a_p:=A x_p$ for every $p \in {P}$ and $c_r:= c^T x_r$ and $a_r := A x_r$ for every $r \in {R}$, we obtain the \textit{aggregated} master problem formulation
\begin{eqnarray}
\min_{\lambda} & \displaystyle \sum_{p \in {P}}c_p\lambda_p + \sum_{r \in {R}}c_r \lambda_r,& \label{eq:agg:obj}\\
\mbox{s.t.} & \displaystyle \ \sum_{p \in {P}}a_p\lambda_p + \sum_{r \in {R}}a_r \lambda_r \leq b,& \label{eq:agg:c1}\\
& \displaystyle \sum_{p \in {P}}\lambda_p = 1,& \label{eq:agg:c2}\\
& \displaystyle \lambda_p \geq 0, & \forall p \in {P}, \label{eq:agg:c3}\\
& \displaystyle \lambda_r \geq 0,  & \forall r \in {R}. \label{eq:agg:c4}
\end{eqnarray}
{Although this is not required, we can still exploit the separability of $\mathcal{X}$ at the subproblem level, without changing \eqref{eq:agg:obj}-\eqref{eq:agg:c4}. In such case, we obtain the extreme points and extreme rays in a decomposed way from each set $\mathcal{X}^k$, $k \in \mathcal{K}$, as in the disaggregated formulation. 
{Then, we use these extreme points or rays together in order to generate a new column. This can be done as any extreme point/ray of $\mathcal{X}$ can be written as a Cartesian product of extreme points/rays from each $\mathcal{X}^k$. 
In order to add this new column, one should use $a_p = A (x^{1}_p, \ldots, x^{K}_p) = A^1 x^{1}_p + \ldots + A^K x^{K}_p$.}
Notice that in the particular case of identical subproblems (i.e. $\mathcal{X}^1 = \ldots = \mathcal{X}^K$), with $A^1 = \ldots = A^K$ and $ c^{1} = \ldots = c^{K}$, we need to solve only one of the subproblems. In such case, it is common to make a variable change in the master variables, as $x^{1}_p = \ldots = x^{K}_p$ and the new column becomes $a_p = A (x^{1}_p, \ldots, x^{1}_p) = K A^1 x^{1}_p$. In the applications that we address in sections \ref{sec:mkl} and \ref{sec:stoch}, namely the multiple kernel learning problem and the two-stage stochastic programming problem, we use aggregated formulations of the master problem as described here.} 

\subsection{Generalized decomposition for convex programming problems}
\label{sub:sec:convex:function}

\begin{figure}[hbt]
\centering
\includegraphics[keepaspectratio=true,scale=0.6,clip=true,viewport = 50 560 540 770]{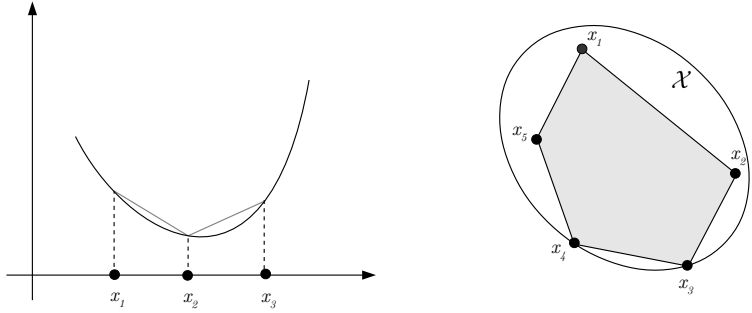}
\caption{Inner linearizations of a convex function and a convex set.}
\label{fig:innerlin}
\end{figure}

As pointed out in \cite{geoffrion1970elements:concepts,geoffrion1970elements:algo} the principles of DWD are based on a more general concept known as \emph{inner linearization}. An inner linearization approximates a convex function by its epigraph so its value at a given point is not underestimated. 
{
Similarly, DWD can be used to approximate a convex set by means of the convex hull of points selected from the set. These points form a \textit{base} so that any point that belongs to the resulting convex hull can be written as a convex combination of the points in the base. 
Hence, we can use the DWD to solve this approximated problem, as observed in \cite[Ch.\ 24]{dantzig1963}. The accuracy of the approximation depends on the points that we include in the base. In Fig.\ref{fig:innerlin} we illustrate inner linearizations for a convex function and a convex set, respectively.
}

\noindent
Let us consider a convex optimization problem defined as follows
\begin{equation}
\label{eq:originalconvex}
\min \ \ f(x), \ \mbox{s.t.} \ \ F(x) \leq 0, \ \ \ x \in \mathcal{X} \subseteq \mathbb{R}^n,
\end{equation}
in which $f : \mathcal{X} \rightarrow \mathbb{R}$ and $F : \mathcal{X} \rightarrow \mathbb{R}^m$ are convex functions which are continuous in $\mathcal{X}$.
For simplicity, we assume that $\mathcal{X}$ is a bounded set and that we have a fine base $(x_1,x_2,\ldots,x_{|Q|}) \in \mathcal{X}$, { where $Q$ is a finite set of indices of points selected from $\mathcal{X}$ \cite{geoffrion1970elements:algo,geoffrion1970elements:concepts}. By fine base, we mean a (finite) set of points that provide an approximation as accurate as we require for $\mathcal{X}$.} Hence, a point $x \in \mathcal{X}$ can be {approximated} by the convex combination of points in $\mathcal{X}$ as
\begin{equation}
\label{eq:rewritexconvex}
x = \sum_{q \in {Q}}  \lambda_q x_q , \mbox{ with } \sum_{q \in {Q}}\lambda^k_q = 1.
\end{equation}
Since the function $f$ in problem (\ref{eq:originalconvex}) is convex, the following relationship must be satisfied for any $x \in \mathcal{X}$
\begin{eqnarray}
f(x) = f \left( \sum_{q \in Q} \lambda_q x_q \right) \leq \sum_{q \in Q} \lambda_q f \left(x_q \right). \label{eq:fxconvexrelationship}
\end{eqnarray}
The right-hand side of \eqref{eq:fxconvexrelationship} is an inner linearization of $f(x)$, which can be used to describe $f(x)$ as closely as desired. As long as we choose the base appropriately, this approximation does not underestimate the value of $f(x)$. The same idea applies to $F(x)$.
By denoting $f_q = f(x_q)$ and $F_q = F(x_q)$, $q \in Q$, we could approximate problem \eqref{eq:originalconvex} as closely as 
require with the following \emph{linear} master problem
\begin{eqnarray}
\min_{\lambda} & \displaystyle \sum_{q \in {Q}} f_q \lambda_q,& \label{eq:convex:obj}\\
\mbox{s.t.} & \displaystyle \ \sum_{q \in {Q}} F_q \lambda_q \leq 0,& \label{eq:convex:c1}\\
& \displaystyle \sum_{q \in {Q}}\lambda_q = 1, & \label{eq:convex:c2}\\
& \displaystyle \lambda_q \geq 0, & \forall q \in {Q}. \label{eq:convex:c3}
\end{eqnarray}
Since the cardinality of $Q$ is typically large, we use the column generation method to solve \eqref{eq:convex:obj}--\eqref{eq:convex:c3}. 

In a given outer iteration, let $\overline{u} \in \mathbb{R}^m_{-}$ and $\overline{v} \in \mathbb{R}$ be the dual optimal solutions associated with constraints (\ref{eq:convex:c1}) and (\ref{eq:convex:c2}), respectively.
We use $\overline{u}$ to call the oracle, which is given by the following convex subproblem
\begin{equation}
{SP}(\overline{u}) := \min_{x \in \mathcal{X}} \ \{ f(x) - \overline{u} F(x) \}.
\end{equation}
{
An optimal solution of this subproblem results in a point of $\mathcal{X}$. 
If ${SP}(\overline{u}) - \overline{v} < 0$, then the point can be added to the current base  
in order to improve the approximation of $\mathcal{X}$. Hence, we obtain a new column of the master problem \eqref{eq:convex:obj}-\eqref{eq:convex:c3}.}
The  multiple kernel learning problem presented in Section \ref{sec:mkl} has a master problem formulation which is similar to (\ref{eq:convex:obj})-(\ref{eq:convex:c3}) and hence we apply the ideas described above.

\section{Stabilized column generation/cutting plane methods} \label{sec:algo}

There is a wide variety of stabilized variants of the column generation/cutting plane methods \cite{Marsten1975,ruszczynsky1986,Kiwiel90,goffin1992,lemarechal1995,bahn1995,
wentges1997,duMerle1999,Frangioni2002,babonneau2007,GonGonMun2013}. 
In this section we briefly present some methodologies which have proven to be very effective for the classes of problems addressed in this paper.
{For the sake of the length of this paper, in this section we only describe a small subset of the many different variants available in the literature.}

\subsection{ACCPM}

The analytic center cutting plane method (ACCPM) proposed in \cite{goffin1992,GofVia02} is an interior point approach that relies on central prices. 
This strategy calculates a dual point which is an approximate analytic center of the localization set associated with the current RMP.
This localization set is given by the intersection of the dual space of the current RMP with the half-space provided by the best lower bound found so far.
Relying on points in the proximity of the analytic center of the localization set usually prevents the unstable behavior between consecutive dual points and also contributes to the generation of fewer and deeper constraints. 
When solving convex optimization problems with nonlinear objective function, the ACCPM can take advantage of second order information to enhance its performance \cite{babonneau2009}. 
An interesting feature of this approach is given by its theoretical 
polynomial
complexity \cite{goffin1996complexity,altmankiwiel96,Kiwiel96}.

\subsection{Bundle methods: Level set}

Bundle methods \cite{Kiwiel90,HirLem1993,Frangioni2002} have become popular techniques to stabilize cutting plane methods. In general, these techniques stabilize the RMP via a proximity control function. There are different classes of bundle methods such as \emph{proximal} \cite{Kiwiel90}, \emph{trust region} \cite{Schramm1992} and \emph{level set} \cite{lemarechal1995}, among others. They use a similar Euclidean-like prox-term to penalize any large deviation in the dual variables. The variants differ in the way they calculate the sequence of iterates. For instance, the level set bundle method is an iterative method which relies on piece-wise linear outer approximations of a convex function. At a given iteration, a level set is created using the best bound found by any of the proposed iterates and the best solution of the outer approximation. Then a linear convex combination is considered to create the level set bound. Finally, the next iterate is obtained by minimizing the prox-term function subject to the structural constraints and the level set bound constraint \cite{lemarechal1995}. 

\subsection{PDCGM}
The primal-dual column generation method (PDCGM) was originally proposed in \cite{GonSar96} and further developed in \cite{GonGonMun2013}. 
A very recent attempt of combining this method in a branch-and-price framework to solve the
vehicle routing problem with time windows can be found in \cite{MunGon2012}.
Warmstarting techniques for this method have been proposed in \cite{Gon98,GonGon12}.
The method relies on sub-optimal and well-centered RMP dual solutions with the aim of reducing the heading-in and tailing-off effects often observed in the standard column generation.
Given a primal-dual feasible solution $(\tilde{\lambda},\tilde{u},\tilde{v})$ of the RMP (\ref{eq:rmp:obj})-(\ref{eq:rmp:c4}), which may have a non-zero distance to optimality,
we can use it to obtain upper and lower bounds for the optimal value of the MP as follows
\begin{eqnarray}
{z}^{ }_{UB}(\tilde{\lambda}) &:=& \sum_{k \in \mathcal{K}} \left( \sum_{p \in \overline{P}^k} c^k_p \tilde{\lambda}^k_p + \sum_{r \in \overline{R}^k} c^k_r \tilde{\lambda}^k_r \right),\\
{z}^{ }_{LB}(\tilde{u},\tilde{v}) &:=& b^{T} \tilde{u} + \sum_{k \in \mathcal{K}} \tilde{v}^k.
\end{eqnarray}
The solution  $(\tilde{\lambda},\tilde{u},\tilde{v})$ is called sub-optimal or $\varepsilon$-optimal solution, if it satisfies
\begin{equation}
0 \leq \left({z}^{ }_{UB}(\tilde{\lambda}) - {z}^{ }_{LB}(\tilde{u},\tilde{v}) \right) \leq \varepsilon ({10^{-10}} + | {z}^{ }_{UB}(\tilde{\lambda})|),
\end{equation}
for some tolerance $\varepsilon > 0$.
Additionally, the primal-dual interior point method provides well-centered dual solutions since it keeps the complementarity products of the primal-dual pairs in the proximity of the central path \cite{Gondzio2011}.
We say a point is well-centered if the following conditions are satisfied
\begin{eqnarray}
\gamma \mu \leq (c^k_p - \tilde{u}^T a^k_p - \tilde{v}^k)\tilde{\lambda}^k_p \leq \frac{1}{\gamma} \mu, \ \forall k \in \mathcal{K}, \ \forall p \in \overline{{P}}^k, \\
\gamma \mu \leq (c^k_r - \tilde{u}^T a^k_r) \tilde{\lambda}^k_r \leq \frac{1}{\gamma} \mu, \ \forall k \in \mathcal{K}, \ \forall r \in \overline{{R}}^k,
\end{eqnarray}
where $\gamma \in (0,1)$ and $\mu$ is the barrier parameter used in the primal-dual interior point algorithm to define the central path \cite{Gondzio2011}.

The PDCGM dynamically adjusts the tolerance used to solve each restricted master problem so it does not stall. The tolerance used to solve the RMPs is loose at the beginning and is tightened as the column generation progresses to optimality.
The method has a very simple description that makes it very attractive (see Algorithm \ref{alg:pdcgm}).

\begin{algorithm}
\normalsize
\caption{The Primal-Dual Column Generation Method}
\label{alg:pdcgm}
\textbf{Input:} Initial RMP; parameters $\kappa$, $\varepsilon_{\max}>0$, $D>1$, $\delta>0$.\\
\textbf{Set:} $\mbox{LB}=-\infty$, $\mbox{UB}=\infty$, $gap = \infty$, $\varepsilon = 0.5$;
\begin{algorithmic}[1]
\WHILE{$( gap \geq \delta )$}
	\STATE{find a well-centered $\varepsilon$-optimal solution $(\tilde{\lambda},\tilde{u},\tilde{v})$ of the RMP;}
	\STATE{$\mbox{UB} = \min\{\mbox{UB}, {z}^{ }_{UB}(\tilde{\lambda})\}$;} \label{ub:store}
	\STATE{call the oracle with the query point $(\tilde{u}, \tilde{v})$;}
	\STATE{$\mbox{LB} = \max\{\mbox{LB}, {z}^{ }_{LB}(\tilde{u}, \tilde{v}) + {z}^{ }_{SP}(\tilde{u}, \tilde{v}) \}$;}
	\STATE{$gap = ( \mbox{UB} - \mbox{LB} ) / ( {10^{-10}} + | \mbox{UB} | )$;}
	\STATE{$\varepsilon = \min \{ \varepsilon_{\max}, \ gap / D \}$;}
	\STATE{if $( {z}^{ }_{SP}(\tilde{u}, \tilde{v}) < 0 )$ then add the new columns to the RMP;}
\ENDWHILE
\end{algorithmic}
\end{algorithm}

\begin{remark}
The PDCGM converges to an optimal solution of the MP in a finite number of outer iterations. 
This has been shown in \cite{GonGonMun2013} based on the valid lower bound provided by an $\varepsilon$-optimal solution and the progress of the algorithm even if columns with nonnegative reduced costs are obtained.
\end{remark}

\begin{remark}
Since close-to-optimality solutions of RMPs are used there is no guarantee of monotonic decrease of the upper bound. Therefore, we update the upper bound using $UB=\min\{UB,{z}^{ }_{UB}(\tilde{\lambda})\}$.
\end{remark}

\begin{remark}
The only parameters PDCGM requires to be set are the \emph{degree of optimality} denoted by $D$ and the initial tolerance threshold $\varepsilon_{max}$.
\end{remark}

Having described the PDCGM and two of the most successful stabilization techniques used in column generation and cutting plane methods, we now describe three different applications and compare the efficiency of these methods on them.

\section{Multiple kernel learning problem (MKL)} \label{sec:mkl}

Kernel based methods are a class of algorithms that are widely used in data analysis. They compute the similarities between two examples $x_j$ and $x_i$ via the so-called kernel function $\kappa(x_j,x_i)$. Typically, the examples are mapped into a feature space by using a mapping function $\Phi$,
so that $\kappa(x_j,x_i) = \langle\Phi(x_j),\Phi(x_i)\rangle$. In practice, one kernel may not be enough to effectively describe the similarities between the data and therefore using multiple kernels provides a more accurate classification. However, this may be more costly in terms of CPU time. Several studies concerning different techniques to solve the multiple kernel learning problem (MKL) are available in the literature, such as \cite{Bach2004,lanckriet2004learning,Sonnenburg2006,Zien2007,rakotomamonjy2008simplemkl,xu2009extended,spicyMKL}. For a thorough taxonomy, classification and comparison of different algorithms developed in this context, see \cite{Gonen2011}. 
In this paper, we are interested in solving a problem equivalent to the semi-infinite linear problem formulation proposed in \cite{Sonnenburg2006} which follows developments in \cite{Bach2004} and that can be solved by a column generation method. We have closely followed the developments of both papers, namely \cite{Bach2004} and \cite{Sonnenburg2006}, keeping a similar notation.
In this section, we describe the single and multiple kernel learning problems in terms of their primal and dual formulations. Then, for the MKL, we derive the column generation components. Finally, we present two computational experiments with the PDCGM for solving publicly available instances and compare these results against several state-of-the-art methods.

\subsection{Problem formulation}

We focus on a particular case of kernel learning known as the support vector machine (SVM) with soft margin loss function \cite{Sonnenburg2006}. As described in \cite{Vapnik98}, the SVM is a classifier proposed for binary classification problems and is based on the theory of structural risk minimization.
The problem can be posed as finding a linear discriminant (hyperplane) with the maximum margin in a given feature space for $n$ data samples $(x_i,y_i)$, where $x_i$ is the $d$-dimensional input vector (features) and $y_i \in \{\pm 1\}$ is the binary class label. We start by describing the single kernel problem and then extend the formulations to the multiple kernel case. In SVM with a single kernel, the discriminant function has usually the following form
\begin{equation}
f_{w,b}(x) = 
\langle w, \Phi(x)\rangle + b , \label{eq:mkl:discriminant}
\end{equation}
where $w \in \mathbb{R}^s$ is the vector of weight coefficients, $b \in \mathbb{R}$ is the bias term and $\Phi: \mathbb{R}^d \rightarrow \mathbb{R}^s$ is the function which maps the examples to the feature space of dimension $s \in \mathbb{Z}_+$.
Hence, we use the sign of $f_{w,b}(x)$ to verify if an example $x$ should be classified as either $-1$ or $+1$.

To obtain the vectors $w$ and $b$ which lead to the best linear discriminant, we can solve the following SVM problem with a single kernel
\begin{eqnarray}
\min_{w,\xi,b} & \displaystyle \frac{1}{2}||w||_2^2 + C\sum_{j \in \mathcal{N}}\xi_j, & \ \ \label{svm:or:obj}\\
\mbox{s.t.} & \quad \displaystyle y_j\left(\langle w, \Phi(x_j)\rangle + b\right) + \xi_j \geq 1, & \ \ \forall j \in \mathcal{N},\label{svm:or:eq1}\\
& \xi_j \geq 0, & \ \ \forall j \in \mathcal{N}, \label{svm:or:eq2}
\end{eqnarray}
where $\mathcal{N} = \{1,\ldots,n\}$, $C \in \mathbb{R}$ is a penalty parameter associated with misclassification and $\xi \in \mathbb{R}^n$ is the vector of variables which measure the error in the classifications. 
{In this formulation, the first term in the objective function aims to maximize the distance between the discriminant function and the features of the training data sample. This distance is called a  \textit{margin} and by maximizing $1/||w||_2^2$ we keep this margin as large as possible \cite{benhur2010}. The second term in the objective function aims to minimize the misclassification of the data sample. The value of parameter $C$ determines the importance of each term for the optimization.}

{Let $\alpha \in \mathbb{R}_+^n$ and $\gamma \in \mathbb{R}_+^n$ denote the vectors of dual variables associated to constraints \eqref{svm:or:eq1} and \eqref{svm:or:eq2}, respectively. To obtain the dual formulation of \eqref{svm:or:obj}--\eqref{svm:or:eq2}, we first define the Lagrangian function:
\begin{equation}
\mathcal{L}(w,\xi,b) = \frac{1}{2}||w||_2^2 + C\sum_{j \in \mathcal{N}}\xi_j - \sum_{j \in \mathcal{N}} \alpha_j \left( y_j\left(\langle w, \Phi(x_j)\rangle + b\right) + \xi_j - 1 \right) - \sum_{j \in \mathcal{N}}\gamma_j \xi_j. \nonumber
\end{equation}

\noindent
{Applying first order optimality conditions to the Lagrangian function and noting that $\alpha \in \mathbb{R}_+^n$ and $\gamma \in \mathbb{R}_+^n$, we obtain the following conditions}
\begin{equation}
w = \sum_{j \in \mathcal{N}} \alpha_j y_j \Phi(x_j), \nonumber
\end{equation}
\begin{equation}
\sum_{j \in \mathcal{N}} \alpha_j y_j = 0, \nonumber
\end{equation}
\begin{equation}
0 \leq \alpha_j \leq C, \ \ \forall j \in \mathcal{N}. \nonumber
\end{equation}
Then, by using these relationships in the Lagrangian function, we obtain the dual problem}
\begin{eqnarray}
\max_{\alpha} & \ \ \ \displaystyle \sum_{j \in \mathcal{N}}\alpha_j - \frac{1}{2}\sum_{j \in \mathcal{N}}\sum_{i \in \mathcal{N}}\alpha_j\alpha_iy_jy_i \kappa(x_j,x_i), & \ \ \label{svm:or_dual:obj}\\
\mbox{s.t.} &  \ \ \ \displaystyle \sum_{j \in \mathcal{N}}\alpha_jy_j = 0, &\label{svm:or_dual:eq1}\\
&  \ \ \ \displaystyle 0 \leq \alpha_j \leq C, & \ \ \forall j \in \mathcal{N}. \label{svm:or_dual:eq2}
\end{eqnarray}

Recall that we use $\kappa(x_j,x_i) = \langle\Phi(x_j),\Phi(x_i)\rangle$ to denote the kernel function which maps $\mathbb{R}^d \times \mathbb{R}^d \rightarrow \mathbb{R}$. This function induces a $n \times n$
\textit{kernel matrix} with entries $\kappa(x_j,x_i)$ for each $i, j \in \mathcal{N}$. Following \cite{Sonnenburg2006}, we only consider kernel functions which lead to positive semi-definite kernel matrices, so \eqref{svm:or_dual:obj}-\eqref{svm:or_dual:eq2} is a convex quadratic programming problem.
Notice that in an optimal solution $\alpha^{*}$ of the dual problem \eqref{svm:or_dual:obj}-\eqref{svm:or_dual:eq2}, $\alpha_j^{*} > 0$ means that the $j$-th constraint in \eqref{svm:or:eq1} is active and, therefore, $x_j$ is a \textit{support vector}.

For many applications and due to the nature of the data, a more flexible approach is to combine different kernels. MKL aims to optimize the kernel weights ($\beta$) while training the SVM. The benefit of using MKL is twofold. On the one hand, it finds the relevant features of a given kernel 
much like the
single kernel learning context, and on the other hand, it leads to an improvement in the classification accuracy since more kernels can be considered. Similar to the single kernel SVM, the discriminant function can be described as
\begin{equation}
f_{w,b,\beta}(x) = 
\sum_{k \in \mathcal{K}} \beta_k \langle w_k, \Phi_k(x)\rangle + b, \label{mkl:discriminant2}
\end{equation}
where $\mathcal{K}$ represents the set of kernels (each corresponding to a positive semi-definite matrix) with a different set of features and $\beta_k$ is the weight associated with kernel $k \in \mathcal{K}$.
We also have that $\beta_k\geq 0$ for every $k \in \mathcal{K}$ and $\sum_{k \in \mathcal{K}}\beta_k = 1$.
Similar to the single kernel learning problem, $w_k$ and $\Phi_k(x)$ are the weight vector and the feature map associated with kernel $k \in \mathcal{K}$, respectively.
The MKL problem was first formulated as a semi-definite programming problem in \cite{lanckriet2004learning}. In \cite{Bach2004}, the authors reformulated this problem as a second-order conic programming problem yielding a formulation which can be solved by sequential minimal optimization (SMO) techniques. 

Similar to the single kernel problem formulation \eqref{svm:or:obj}-\eqref{svm:or:eq2}, the MKL primal problem for classification can be formulated as 
\begin{eqnarray}
\min_{v,\xi,b} & \ \ \ \displaystyle \frac{1}{2}\left(\sum_{k \in \mathcal{K}}||v_k||_2\right)^2 + C\sum_{j \in \mathcal{N}}\xi_j, & \ \ \label{mkl:or:obj}\\
\mbox{s.t.} &  \ \ \ \displaystyle y_j\left(\sum_{k \in \mathcal{K}}\langle v_k, \Phi_k(x_j)\rangle + b\right) + \xi_j \geq 1, & \ \ \forall j \in \mathcal{N},\label{mkl:or:eq1}\\
&  \ \ \ \displaystyle  \xi_j \geq 0, & \ \ \forall j \in \mathcal{N}, \label{mkl:or:eq2}
\end{eqnarray}
where $v_k \in \mathbb{R}^{s_{k}}$ 
and $s_{k}$ is the dimension
of the feature space associated with kernel $k$ and $v_k = \beta_k w_k$, for every $k \in \mathcal{K}$. As pointed out in \cite{Zien2007}, the use of $v$ instead of $\beta$ and $w$ (as in \eqref{mkl:discriminant2}) makes problem \eqref{mkl:or:obj}-\eqref{mkl:or:eq2} convex.
This problem can be interpreted as (i) finding a linear convex combination of kernels while (ii) maximizing the normalized distance between the features of the training data sample and the discriminant function and (iii) minimizing the misclassification of the data sample for each kernel function. 
Following \cite{Bach2004}, we associate $\alpha  \in \mathbb{R}_+^n$ with constraint \eqref{mkl:or:eq1} so the dual formulation of \eqref{mkl:or:obj}-\eqref{mkl:or:eq2} can be written as 
\begin{eqnarray}
\max_{\alpha,\rho} & \ \ \ \displaystyle - \rho, & \label{mkl:dual1:obj}\\
\mbox{s.t.} & \ \ \ \displaystyle S^k(\alpha) - \rho \leq 0, &\ \ \forall k \in \mathcal{K}, \label{mkl:dual1:eq1}\\
& \ \ \ \displaystyle \alpha \in \Gamma, &\label{mkl:dual1:eq2}
\end{eqnarray}
where for every $k \in \mathcal{K}$
\begin{equation}
S^k(\alpha):= \frac{1}{2}\sum_{j \in \mathcal{N}} \sum_{i\in \mathcal{N}}\alpha_j\alpha_iy_jy_i\kappa_k(x_j,x_i)-\sum_{j \in \mathcal{N}}\alpha_j, \label{eq:skalpha}
\end{equation}
$\displaystyle \kappa_k(x_j,x_i) = \langle\Phi_k(x_j),\Phi_k(x_i)\rangle$ and 
\begin{equation}
\Gamma = \left\{\alpha \mid \sum_{j \in \mathcal{N}} \alpha_j y_j = 0, 0 \leq \alpha_j \leq C, \forall j \in \mathcal{N}\right\}.
\end{equation}
Note that problem \eqref{mkl:dual1:obj}-\eqref{mkl:dual1:eq2} is a quadratically constrained quadratic problem (QCQP) \cite{loboSOCP1998}.
Since the number of training examples and kernel matrices used in practice are typically large, it may become a very challenging large-scale problem. Nevertheless, this problem can be effectively solved if we use the inner linearization technique addressed in Section \ref{sub:sec:convex:function}.

\subsection{Decomposition and column generation formulation}

In this section we derive the master problem formulation and the oracle which are associated with problem \eqref{mkl:dual1:obj}-\eqref{mkl:dual1:eq2}. All the developments closely follow Section \ref{sec:dwdcolgen}.
Let $P$ be the set of indices of points in the interior and boundaries of 
the
set $\Gamma$. 
Since $\Gamma$ is a bounded set, we can write any $\alpha \in \Gamma$ as a convex combination of the points $\alpha_p$, $p \in {P}$.
We need extreme as well as interior points of 
$\Gamma$ as $S_k(\alpha)$ is a convex quadratic function and therefore an attractive $\alpha$ may lie in the interior of 
$\Gamma$.
Using the developments presented in Section \ref{sub:sec:convex:function} the following master problem is associated with \eqref{mkl:dual1:obj}-\eqref{mkl:dual1:eq2}
\begin{eqnarray}
\min_{\rho,\lambda} & \ \ \ \rho, & \label{mkl:mp:obj}\\
\mbox{s.t.} & \ \ \ \displaystyle 
\rho - \sum_{p \in {P}}\lambda_pS^k(\alpha_p) \geq 0, & \ \ \forall k \in \mathcal{K},\label{mkl:mp:eq1}\\
& \ \ \ \displaystyle \sum_{p \in {P}} \lambda_p = 1, &\label{mkl:mp:eq2}\\
&  \ \ \ \displaystyle \lambda_p \geq 0, &\ \ \forall p \in {P}.\label{mkl:mp:eq3}
\end{eqnarray}
Note that we have changed the direction of the optimization problem and now we minimize $\rho$ (compare \eqref{mkl:dual1:obj} with \eqref{mkl:mp:obj}). This change does not affect the column generation, but the value of the objective function will have the opposite sign. Also, observe that $\rho \in \mathbb{R}$ only appears in the master problem.
It is worth mentioning that this master problem formulation is the dual of the semi-infinite linear problem presented in \cite{Sonnenburg2006}. 
Due to the large number of possible points $\alpha_p \in \Gamma$, we rely on column generation to solve the MKL.
The corresponding oracle is given by a single pricing subproblem. Let $\overline{\beta} \in \mathbb{R}^{|\mathcal{K}|}_{+}$ and $\overline{\theta} \in \mathbb{R}$ be the dual solutions associated with constraints \eqref{mkl:mp:eq1} and \eqref{mkl:mp:eq2}, respectively.
We always have $\overline{\beta}_k \geq 0$ for all $k \in \mathcal{K}$. In addition, from the dual of problem \eqref{mkl:mp:obj}-\eqref{mkl:mp:eq3}, we have that $\sum_{k \in \mathcal{K}}\overline{\beta}_k = 1$. 
Indeed, it can be shown that the duals $\overline{\beta}_k$ are the weights associated with the kernels in $\mathcal{K}$ \cite{Sonnenburg2006}. 
The subproblem is defined as
\begin{equation}
{
{SP}(\overline{\beta}) := \min_{\alpha \in \Gamma} \left\{\sum_{k \in \mathcal{K}} \overline{\beta}_kS^k(\alpha)\right\}.}\label{mkl:subproblem}
\end{equation}
This subproblem turns 
out
to be of the same form as a single kernel problem such as  \eqref{svm:or_dual:obj}-\eqref{svm:or_dual:eq2}, with $\kappa(x_j,x_i) = \sum_{k \in \mathcal{K}}\beta_k \kappa_k(x_j,x_i)$. Therefore, an SVM solver can be used to solve the subproblem.
Finally, the value of the oracle is $z_{SP}(\overline{\beta},\overline{\theta}) := \min \{0;SP(\overline{\beta})-\overline{\theta}\},$
and a new column is obtained if $z_{SP}(\overline{\beta},\overline{\theta}) <0$.

\subsection{Computational experiments}

To evaluate the efficiency of the PDCGM for solving the MKL, we have carried out computational experiments using benchmarking instances from the UCI Machine Learning Repository data sets \cite{UCI}.
The pricing subproblem \eqref{mkl:subproblem} is solved by a state-of-the-art machine learning toolbox, namely SHOGUN 2.0.0 \cite{SHOGUN2010} (\url{http://www.shogun-toolbox.org/}). The first set of experiments replicates the settings proposed in \cite{rakotomamonjy2008simplemkl}. Hence, we have selected 3 polynomial kernels of degrees 1 to 3 and 10 Gaussian kernels with 
widths $2^{-3},2^{-2},\ldots,2^{6}$.
We have used these 13 kernels to generate $13(d+1)$ kernel matrices by using all features and every single feature separately, where $d$ is the number of features of the instance. 
We have selected $70\%$ of the data as 
a
training sample and normalised it to unit trace. The remaining data (\emph{i.e.,} $30\%$) is used for testing the accuracy of the approach. We have randomly generated 20 instances per data set. We have chosen $C = 100$ as penalty parameter and the algorithm stops when the relative duality gap drops below $\delta = 10^{-2}$ (standard stopping criterion for this application).
The degree of optimality was set to $D = 10$ and the initial RMP had only one column which was generated by solving problem \eqref{mkl:subproblem} with equal weights (\emph{i.e.,} $\overline{\beta}_k = \frac{1}{|\mathcal{K}|}$).
A Linux PC with an Intel Core i5 2.3 GHz CPU and 4.0 GB of memory 
was used to run the experiments {(single thread)}. 

In Table \ref{table:mkl:exp1} we compare the performance of the PDCGM with the performances of 
three other methods presented in \cite{rakotomamonjy2008simplemkl}, namely: 
\begin{itemize}
 \item SILP: The Semi-infinite linear programming algorithm was proposed in \cite{Sonnenburg2006}; {therein} a cutting plane method without stabilization {is used to solve}
 the dual of \eqref{mkl:mp:obj}-\eqref{mkl:mp:eq3}.
 \item simpleMKL: The simple subgradient descent method was proposed in \cite{rakotomamonjy2008simplemkl}. The method solves problem \eqref{mkl:or:obj}-\eqref{mkl:or:eq2} with a weighted $l_2$-norm regularization and an extra constraint for the kernel weights. This constrained optimization problem is solved by a gradient method on the simplex. The simpleMKL updates the descent direction by looking for the maximum admissible step and only calculates the gradient of the objective function when the objective function stops decreasing.
\item GD: This is the standard gradient descent method which solves the same constrained optimization problem as the simpleMKL. In this approach, a descent direction is recomputed every time the kernel weights are updated.
\end{itemize}

\begin{table}[tbp]

\setlength{\tabcolsep}{0.7ex}
\caption{MKL: comparison between PDCGM results and results reported in \cite{rakotomamonjy2008simplemkl} for the SILP, simpleMKL and GD when using 70\% of the data as training sample and 20 instances per data set.}
\scriptsize
\centering
    \begin{tabular}{cccccc}
\hline \noalign{\smallskip}
\multicolumn{1}{c}{{Instance}} &	\multicolumn{1}{c}{Method} & \multicolumn{1}{c}{{\# Kernel}} & \multicolumn{1}{c}{{Accuracy(\%)}} & \multicolumn{1}{c}{{Time(s)}} & \multicolumn{1}{c}{\# SVM calls}\\
\hline \noalign{\smallskip} 
{bupa/liver}			&	PDCGM	&	11.6	$\pm$	2.2	&	68.6	$\pm$	3.0	&	6.7	$\pm$	0.6	&	10.1	$\pm$	0.7	\\
$n	=	242$	&	SILP	&	10.6	$\pm$	1.3	&	65.9	$\pm$	2.6	&	23.8	$\pm$	4.9	&	99.8	$\pm$	20.0	\\
$|\mathcal{K}|	=	91$	&	GD	&	11.6	$\pm$	1.3	&	66.1	$\pm$	2.7	&	15.7	$\pm$	7.1	&	972.0	$\pm$	630.0	\\
			&	simpleMKL	&	11.2	$\pm$	1.2	&	65.9	$\pm$	2.3	&	9.5	$\pm$	6.3	&	522.0	$\pm$	382.0	\\
\hline	\noalign{\smallskip}																				
{ionosphere}			&	PDCGM	&	10.2	$\pm$	0.9	&	95.3	$\pm$	2.0	&	41.4	$\pm$	4.8	&	16.6	$\pm$	1.8	\\
$n	=	246$	&	SILP	&	21.6	$\pm$	2.2	&	91.7	$\pm$	2.5	&	267.5	$\pm$	52.5	&	95.6	$\pm$	13.0	\\
$|\mathcal{K}|	=	442$	&	GD	&	22.9	$\pm$	3.2	&	92.1	$\pm$	2.5	&	210.5	$\pm$	31.0	&	873.0	$\pm$	147.0	\\
			&	simpleMKL	&	23.6	$\pm$	2.6	&	91.5	$\pm$	2.5	&	61.5	$\pm$	23.0	&	314.0	$\pm$	44.0	\\
\hline	\noalign{\smallskip}																				
{pima}			&	PDCGM	&	9.2	$\pm$	3.0	&	76.5	$\pm$	1.6	&	41.4	$\pm$	4.1	&	9.6	$\pm$	0.9	\\
$n	=	538$	&	SILP	&	11.6	$\pm$	1.0	&	76.5	$\pm$	2.3	&	112	$\pm$	18.5	&	157.0	$\pm$	44.0	\\
$|\mathcal{K}|	=	117$	&	GD	&	14.8	$\pm$	1.4	&	75.5	$\pm$	2.5	&	109.5	$\pm$	12.0	&	2620.0	$\pm$	232.0	\\
			&	simpleMKL	&	14.7	$\pm$	1.4	&	76.5	$\pm$	2.6	&	39.5	$\pm$	6.5	&	618.0	$\pm$	148.0	\\
\hline	\noalign{\smallskip}																				
{sonar}			&	PDCGM	&	9.8	$\pm$	2.4	&	88.1	$\pm$	2.9	&	22.6	$\pm$	2.9	&	10.6	$\pm$	1.5	\\
$n	=	146$	&	SILP	&	33.5	$\pm$	3.8	&	80.5	$\pm$	5.1	&	1145.0	$\pm$	432.0	&	403.0	$\pm$	53.0	\\
$|\mathcal{K}|	=	793$	&	GD	&	35.7	$\pm$	3.9	&	80.2	$\pm$	4.7	&	234.5	$\pm$	45.0	&	4000.0	$\pm$	874.0	\\
			&	simpleMKL	&	36.7	$\pm$	5.1	&	80.6	$\pm$	5.1	&	81.5	$\pm$	46.5	&	1170.0	$\pm$	369.0	\\
\hline	\noalign{\smallskip}																				
{wpbc}			&	PDCGM	&	24.0	$\pm$	16.8	&	80.0	$\pm$	4.1	&	9.7	$\pm$	1.0	&	8.9	$\pm$	0.9	\\
$n	=	138$	&	SILP	&	13.7	$\pm$	2.5	&	76.8	$\pm$	1.2	&	44.3	$\pm$	16.0	&	903.0	$\pm$	187.0	\\
$|\mathcal{K}|	=	442$	&	GD	&	16.8	$\pm$	2.8	&	76.9	$\pm$	1.5	&	53.0	$\pm$	3.1	&	7630.0	$\pm$	2600.0	\\
			&	simpleMKL	&	15.8	$\pm$	2.4	&	76.7	$\pm$	1.2	&	10.3	$\pm$	3.1	&	2770.0	$\pm$	1560.0	\\
\hline \noalign{\smallskip}
    \end{tabular}
  \label{table:mkl:exp1}
\end{table}

The first column in Table \ref{table:mkl:exp1} provides the name, size of the training sample ($n$) and number of kernels $(|\mathcal{K}|)$ for each class of instances. For each method described in the second column, we show the average results over 20 randomly generated instances of the corresponding class. The first number represents the average while the second, the standard deviation. Columns 3 to 6 show the number of kernels with no vanishing values $\beta$ at the optimal solution (\# Kernel), the accuracy (Accuracy(\%)), the CPU time in seconds (Time(s)) and the number of calls to the support vector machine solver (\# SVM calls). Accuracy corresponds to the percentage of correct classifications made by the resulting discriminant function. 

From the results in Table \ref{table:mkl:exp1}, we can observe that the PDCGM solves all the instances in less than $42$ seconds on average. Also, it shows a variability of around $10\%$ for the average time. Below the PDCGM results, we have included the results presented in \cite{rakotomamonjy2008simplemkl} for the same statistics. {For columns \textit{\#Kernel}, \textit{Accuracy} and \textit{\#SVM calls} we have taken the results exactly as they are given in \cite{rakotomamonjy2008simplemkl}, while the values in the column \textit{Time} have been scaled by a factor of $0.5$ with respect to the originals. This is because these experiments were run on a Intel Pentium D 3 GHz CPU and 3 GB of memory. Considering the benchmark available at \url{https://www.cpubenchmark.net/singleThread.html}, this machine has a score of 695 whereas the machine we used has a score of 1331. Nevertheless, these CPU times are provided for information only and indicate that PDCGM is competitive regarding CPU times with respect to the other methods. }
In addition to that, the results reported in the last column of Table \ref{table:mkl:exp1} clearly demonstrate that the PDCGM requires fewer calls of the SVM solver than any of the other three methods. Particularly, if we compare the SILP and the PDCGM results, we can observe how much can be gained by using a natural stabilized column generation method over the unstabilized version. Also, {the method seems to be at least as effective as the other methods since it provides a similar or higher level of accuracy and a smaller number of kernels for most of the data sets. This indicates that PDCGM could find more efficient weights to combine the available kernels.}
Finally, an important observation is that the simpleMKL can warmstart the SVM solver due to the small variation from iteration to iteration in the kernel weights and therefore an excessively large number of calls to the SVM solver does not translate into an excessively large CPU time \cite{rakotomamonjy2008simplemkl}. The other three methods (SILP, GD and PDCGM) do not rely on this feature.

In a second set of computational experiments, we have used the same database \cite{UCI} to generate instances following the experiments presented in \cite{xu2009extended}. The only difference with the previous experiments is that we now use $50\%$ of the data available as training sample and the other $50\%$ for testing purposes. 
This second set of experiments allows us to compare our approach against the SILP and simpleMKL (already described) and the extended level method (LSM), a state-of-the-art implementation to solve the MKL{\cite{xu2009extended}}. The LSM belongs to a family of bundle methods and it considers information of the gradients of all previous iterations (as the SILP) and computes the projection onto a level set as a way of regularizing the problem (as the simpleMKL). The results of the PDCGM and the results reported in \cite{xu2009extended} solving the MKL problem are presented in Table \ref{table:mkl:exp2}.

\begin{table}[tbp]
\setlength{\tabcolsep}{2.5ex}
\caption{MKL: comparison between PDCGM results and results reported in \cite{xu2009extended} for the SILP, simpleMKL and LSM when using 50\% of the data as training sample and 20 instances per data set.}
\scriptsize
\centering
\begin{tabular}{ccccc}
\hline \noalign{\smallskip}
 \multicolumn{1}{c}{{Instances}} & \multicolumn{1}{c}{Method} & \multicolumn{1}{c}{{\# Kernel}} & \multicolumn{1}{c}{{Accuracy(\%)}} & \multicolumn{1}{c}{{Time(s)} }\\
\hline \noalign{\smallskip}
{breast} & PDCGM & 8.1 $\pm$ 1.4 & 96.9 $\pm$ 0.7 & 4.2 $\pm$ 0.5 \\
$n = 342$  & SILP  & 10.6 $\pm$ 1.1 & 96.6 $\pm$ 0.8 & 54.2 $\pm$ 9.4\\
$|\mathcal{K}| = 130$ & simpleMKL & 13.1 $\pm$ 1.7 & 96.6 $\pm$ 0.9 & 47.4 $\pm$ 8.9\\
& LSM & 13.3 $\pm$ 1.5 & 96.6 $\pm$ 0.8 & 4.6 $\pm$ 1.0\\
\hline \noalign{\smallskip}
{heart} & PDCGM & 10.2 $\pm$ 3.8 & 82.4 $\pm$ 2.1 & 3.5 $\pm$ 0.2\\
$n = 135$ & SILP  & 15.2 $\pm$ 1.5 & 82.2 $\pm$ 2.0 & 79.2 $\pm$ 38.1\\
$|\mathcal{K}| = 182$ & simpleMKL    & 17.5 $\pm$ 1.8 & 82.2 $\pm$ 2.2 & 4.7 $\pm$ 2.8\\
 & LSM & 18.6 $\pm$ 1.9 & 82.2 $\pm$ 2.1 & 2.1 $\pm$ 0.4\\
\hline \noalign{\smallskip}
{ionosphere} & PDCGM & 12.8 $\pm$ 9.2 & 95.0 $\pm$ 1.5 & 20.9 $\pm$ 3.6\\
$n = 176$  & SILP  & 24.4 $\pm$ 3.4 & 92.0 $\pm$ 1.9 & 1161.0 $\pm$ 344.2\\
$|\mathcal{K}| = 455$ & simpleMKL & 26.9 $\pm$ 4.0 & 92.1 $\pm$ 2.0 & 33.5 $\pm$ 11.6\\
    & LSM & 25.4 $\pm$ 3.9 & 92.1 $\pm$ 1.9 & 7.1 $\pm$ 4.3\\
\hline \noalign{\smallskip}
{pima} & PDCGM & 10.1 $\pm$ 3.5 & 76.4 $\pm$ 0.9 & 20.1 $\pm$ 1.1\\
$n = 384$ & SILP  & 12.0 $\pm$ 1.8 & 76.9 $\pm$ 2.1 & 62.0 $\pm$ 15.2\\
$|\mathcal{K}| = 117$ & simpleMKL & 16.6 $\pm$ 2.2 & 76.9 $\pm$ 1.9 & 39.4 $\pm$ 8.8\\
    & LSM & 17.6 $\pm$ 2.6 & 76.9 $\pm$ 2.1 & 9.1 $\pm$ 1.6 \\
\hline \noalign{\smallskip}
{sonar} & PDCGM & 9.4 $\pm$ 2.8 & 86.4 $\pm$ 2.5 & 11.1 $\pm$ 2.3\\
$n = 104$ & SILP  & 34.2 $\pm$ 2.6 & 79.3 $\pm$ 4.2 & 1964.3 $\pm$ 68.4\\
$|\mathcal{K}| = 793$ & simpleMKL    & 39.8 $\pm$ 3.9 & 79.1 $\pm$ 4.5 & 60.1 $\pm$ 29.6\\
 & LSM & 38.6 $\pm$ 4.1 & 79.0 $\pm$ 4.7 & 24.9 $\pm$ 10.6\\
\hline \noalign{\smallskip}
{vote} &	PDCGM & 10.4 $\pm$ 1.2 & 95.5 $\pm$ 0.9 & 4.9 $\pm$ 0.4\\
$n = 218$ & SILP  & 10.6 $\pm$ 2.6 & 95.7 $\pm$ 1.0 & 26.3 $\pm$ 12.4\\
$|\mathcal{K}| = 221$ & simpleMKL & 14.0 $\pm$ 3.6 & 95.7 $\pm$ 1.0 & 23.7 $\pm$ 9.7\\
	& LSM & 13.8 $\pm$ 2.6 & 95.7 $\pm$ 1.0 & 4.1 $\pm$ 1.3\\
\hline \noalign{\smallskip}
{wdbc} & PDCGM & 8.9 $\pm$ 1.8 & 96.8 $\pm$ 0.9 & 13.1 $\pm$ 1.6\\
$n = 285$ & SILP  & 12.9 $\pm$ 2.3 & 96.5 $\pm$ 0.9 & 146.3 $\pm$ 48.3\\
$|\mathcal{K}| = 403$ & simpleMKL & 16.6 $\pm$ 3.2 & 96.7 $\pm$ 0.8 & 122.9 $\pm$ 38.2\\
	& LSM & 15.6 $\pm$ 3.0 & 96.7 $\pm$ 0.8 & 15.5 $\pm$ 7.5\\
\hline \noalign{\smallskip}
{wpbc} & PDCGM & 9.1 $\pm$ 1.3 & 77.7 $\pm$ 2.8 & 4.2 $\pm$ 0.5 \\
$n = 99$ & SILP  & 17.2 $\pm$ 2.2 & 76.9 $\pm$ 2.8 & 142.0 $\pm$ 122.3\\
$|\mathcal{K}| = 442$ & simpleMKL    & 19.5 $\pm$ 2.8 & 77.0 $\pm$ 2.9 & 7.8 $\pm$ 2.4\\
    & LSM & 20.3 $\pm$ 2.6 & 76.9 $\pm$ 2.9 & 5.3 $\pm$ 1.3 \\
\hline \noalign{\smallskip}
    \end{tabular}
  \label{table:mkl:exp2}
\end{table}
In this case we have omitted the number of calls to the SVM solver since this information was not reported in \cite{xu2009extended}. 
Comparing Table \ref{table:mkl:exp1} and Table \ref{table:mkl:exp2} one can observe that for the same databases, namely \texttt{ionosphere}, \texttt{pima}, \texttt{sonar} and \texttt{wpbc}, the time and accuracy have decreased when considering $50\%$ of training sample instead of $70\%$. These results are expected since the problem sizes are smaller when a $50\%$ training sample is used. 
The PDCGM solves all the instances in less than $21$ seconds on average showing again a small variability on this metric.
{
The machine used in \cite{xu2009extended} has similar characteristics to the one used to run our experiments (Intel 3.2 GHz and 2 GB of RAM memory). Therefore, we have kept the results as they appeared in \cite{xu2009extended}
}

The results obtained with the PDCGM demonstrate that the method is reliable with a performance similar to the state-of-the-art solver, such as the LSM.
Moreover, the performance of the PDCGM seems to be more stable than those of the other methods showing the smallest standard deviation on CPU times for every instance set considered.
Similar to the previous results, and in terms of number of kernels with non-vanishing weights and accuracy, the PDCGM seems to be as effective as the SILP, simpleMKL and the LSM. 

In additional experiments, we have tested the performance of the PDCGM and the LSM (bundle method), leaving aside the influence of the SVM solvers. 
The LSM implementation uses the optimization toolbox MOSEK (\url{http://www.mosek.com}) to solve the auxiliary subproblems. 
As the results indicate, the PDCGM shows less variability than the LSM regarding the CPU time spent on solving the RMPs.
In addition, the PDCGM bottleneck is in the SVM solver (solving the subproblem) unlike the LSM where the bottleneck is in building the cutting plane model and computing the projection onto the level set (solving the restricted master problem). Indeed, the experiments show that the time required by the PDCGM to solve the RMPs does not exceed $2\%$ of the total CPU time on average. Therefore, combining the PDCGM with a more efficient and well-tuned SVM solver implementation may lead to large savings with respect to the state-of-the-art approaches.
\section{Two-stage stochastic programming problem (TSSP)}\label{sec:stoch}

The stochastic programming field has grown steadily in the last decades. Currently, stochastic programming problems can be used to formulate a wide range of applications which arise in real-life situations. 
For an introduction to stochastic programming we refer the interested reader to \cite{Birge1997,Kall1994}.
The two-stage stochastic linear programming problem (TSSP) deals with uncertainty through the analysis of possible outcomes (scenarios).
Several solution methodologies based on column generation and cutting plane methods have been proposed to deal with this class of problems \cite{vanslyke1969,ruszczynsky1986,birge1988,bahn1995,zverovichetal2012}. 
The TSSP can be posed as two interconnected problems. One of them called the \textit{recourse} problem containing only the \textit{second-stage} variables and the other one called the \textit{first-stage} problem which takes into consideration information related to the recourse problem but only optimizes over the set of first stage variables.

The first-stage problem can be stated as 
\begin{eqnarray}
\min_{x} & \displaystyle c^Tx + E[Q(x,\omega)],& \ \ \label{stoch:prog:obj}\\
\mbox{s.t.} & \displaystyle Ax  = b, &\label{stoch:prog:eq1}\\
& \displaystyle x \geq 0, &\label{stoch:prog:eq2}
\end{eqnarray}
where $c \in \mathbb{R}^n$ is the vector of cost coefficients, $x \in \mathbb{R}^n$ is the vector of first-stage variables (which are scenario independent), $\omega \in \Omega$ is the realization of a random event and $\Omega$ the set of possible events, $A \in \mathbb{R}^{m \times n}$ is a scenario independent matrix and $b \in \mathbb{R}^n$ the corresponding right hand side. 
Additionally, $E[Q(x,\omega)]$ represents the expected cost of all possible outcomes, while $Q(x,\omega)$ is the optimal value of the second-stage problem defined by
\begin{eqnarray}
\min_{y} & \ \ \displaystyle q(\omega)^Ty(\omega), & \ \ \label{stoch:prog:sec:obj}\\
\mbox{s.t.} & \ \ \displaystyle W(\omega)y(\omega) = h(\omega) - T(\omega)x, &\label{stoch:prog:sec:eq1}\\
& \ \ \displaystyle y(\omega) \geq 0, &\label{stoch:prog:sec:eq2}
\end{eqnarray}
where $y$ is the vector of second-stage variables and $W$, $T$, $h$ and $q$ are matrices and vectors which depend on the realization of the random event $\omega$. 

\subsection{Problem formulation}
If we assume that set $\Omega$ is discrete and that every realization has a probability of occurrence associated with it, the TSSP problem can be stated as a deterministic equivalent problem (DEP) of the following form 
\begin{eqnarray}
\min_{x,y} & \displaystyle c^Tx + \sum_{i \in \mathcal{S}} p_i q_i^Ty_i, & \ \ \label{stoch:primal:obj}\\
\mbox{s.t.} & \displaystyle Ax  = b, &\label{stoch:primal:eq1}\\
& \displaystyle T_ix + W_i y_i = h_i, & \forall i \in \mathcal{S},\label{stoch:primal:eq2}\\
& \displaystyle x \geq 0, &\label{stoch:primal:eq3}\\
& \displaystyle y_i \geq 0, & \forall i \in \mathcal{S},\label{stoch:primal:eq4}
\end{eqnarray}
where $\mathcal{S}$ is the set of possible scenarios and $p_i$ is the probability that scenario $i$ occurs, with $p_i > 0$ for every scenario $i \in \mathcal{S}$. Also, $q_i \in \mathbb{R}^{\tilde{n}}$ is the column vector cost associated with the second-stage variables $y_i \in \mathbb{R}^{\tilde{n}}$ for every $i \in \mathcal{S}$. For completeness we also define $T_i \in \mathbb{R}^{\tilde{m} \times n}$, $W_i \in \mathbb{R}^{\tilde{m} \times \tilde{n}}$ and $h_i \in \mathbb{R}^{\tilde{m}}$ for every $i \in \mathcal{S}$.

When modeling real-life situations, the DEP is a large-scale problem.
Due to the structure of $A$, $T_i$ and $W_i$, this formulation has an L-shaped form \cite{vanslyke1969} and can be solved using Benders decomposition \cite{Ben62}, 
which can be seen as the dual counterpart of the DWD. Thus, 
to obtain a suitable structure for applying the DWD, we need to work on the dual of problem \eqref{stoch:primal:obj}-\eqref{stoch:primal:eq4}. By associating the vector of dual variables $\eta \in \mathbb{R}^m$ to constraints \eqref{stoch:primal:eq1} and $\theta_i \in \mathbb{R}^{\tilde{m}}$ to \eqref{stoch:primal:eq2}, for every $i \in \mathcal{S}$, the dual of the DEP can be stated as
\begin{eqnarray}
\max_{\eta,\theta} & \displaystyle b^T\eta + \sum_{i \in \mathcal{S}} h_i^T\theta_i, & \ \ \label{stoch:dual:obj}\\
\mbox{s.t.} & \displaystyle A^T\eta  + \sum_{i \in \mathcal{S}} T_i^T \theta_i \leq c, &\label{stoch:dual:eq1}\\
& \displaystyle W_i^T\theta_i \leq p_iq_i, & \forall i \in \mathcal{S}.\label{stoch:dual:eq2}
\end{eqnarray}
Note that in both problems, the primal and the dual, the number of constraints is 
proportional to the number of scenarios and therefore decomposition techniques are often used to solve them. The number of constraints of the dual (primal) problem is $n + \tilde{n}|\mathcal{S}|$ ($m + \tilde{m}|\mathcal{S}|$). One may attempt to solve 
any of these problems
directly with a simplex-type or interior point method; however, it has been shown that by using decomposition approaches one can solve very large problems whereas a direct method may fail due to the number of scenarios considered \cite{vanslyke1969,zverovichetal2012}.

\subsection{Decomposition and column generation formulation}

The master problem associated with DEP formulation
\eqref{stoch:dual:obj}-\eqref{stoch:dual:eq2} is derived in this section. Also, we describe the subproblem which turns 
out
to be the dual of the recourse problem \eqref{stoch:prog:sec:obj}-\eqref{stoch:prog:sec:eq2} and decomposable by scenario. As in \cite{zverovichetal2012}, we are interested in the aggregated version of the master problem. 

The problem decomposes naturally in scenarios having \eqref{stoch:dual:eq1} as the linking constraints. Let us define 
\begin{equation}
\Theta^i = \left\{\theta_i \mid W_i^T\theta_i \leq q_i\right\}, \forall i \in \mathcal{S}, \label{theta:set}
\end{equation}
and $\Theta = \Theta^1\times\Theta^2\times\ldots\times\Theta^{|\mathcal{S}|}$. Note the slight difference between \eqref{theta:set} and the set defined by constraints \eqref{stoch:dual:eq2}. This is similar to the substitution used in \cite{dantzigmadansky1961} where the dual variable associated with the constraint in \eqref{theta:set} for a given scenario is divided by their corresponding probability, $p_i$. We assume that \eqref{theta:set} describes a non-empty convex polyhedron. Note that unlike the other two problems studied in this paper (\emph{i.e.,} MKL and MCNF), the set defined by \eqref{theta:set} is not bounded and therefore we also need extreme rays to write an equivalent Dantzig-Wolfe reformulation. 

Let $\mathcal{P}$ and $\mathcal{R}$ denote the index sets of extreme points and extreme rays of $\Theta$, respectively. We can write any point $\theta \in \Theta$ in terms of these extreme points and rays. Since we exploit the separability of $\Theta$, we have $\theta = (\theta_1, \theta_2, \ldots, \theta_{|\mathcal{S}|})$, with $\theta_i \in \Theta^i$, $i \in \mathcal{S}$. Hence, any extreme point (ray) of $\Theta$ corresponds to the Cartesian product of extreme points (rays) of each $\Theta^i$, $i \in \mathcal{S}$, so we have
\begin{align}
&\theta_i:= \sum_{p \in \mathcal{P}}\lambda_p\left(p_i \theta_p^i\right) + \sum_{r \in \mathcal{R}}\lambda_r\left(p_i \theta_r^i\right), \ \mbox{ with } \ \displaystyle \sum_{p \in \mathcal{P}}\lambda_p = 1, \nonumber
\end{align}
where $\lambda_p \geq 0, \forall p \in \mathcal{P}, \lambda_r \geq 0, \forall r \in \mathcal{R}$, and $\theta_p^i$ and $\theta_r^i$  {correspond to} extreme points and rays of set $\Theta^i$, $i \in \mathcal{S}$, respectively. 
Therefore, the aggregated master problem associated with \eqref{stoch:dual:obj}-\eqref{stoch:dual:eq2} is
\begin{eqnarray}
\max_{\eta,\lambda} & \displaystyle b^T\eta + \sum_{i \in \mathcal{S}}h_i^T\left(\sum_{p \in \mathcal{P}}\lambda_p\left(p_i \theta_p^i\right) + \sum_{r \in \mathcal{R}}\lambda_r\left(p_i \theta_r^i\right)\right), & \label{stoch:mp:agg:obj}\\
\mbox{s.t.} & \displaystyle A^T\eta + \sum_{i \in \mathcal{S}}T_i^T\left(\sum_{p \in \mathcal{P}}\lambda_p\left(p_i\theta_p^i\right) + \sum_{r \in \mathcal{R}}\lambda_r\left(p_i\theta_r^i\right)\right) \leq c, &\label{stoch:mp:agg:eq1}\\
& \displaystyle \sum_{p \in \mathcal{P}} \lambda_p = 1, &\label{stoch:mp:agg:eq2}\\
&  \displaystyle \lambda_p \geq 0, & \! \! \! \! \forall p \in \mathcal{P},\label{stoch:mp:agg:eq3}\\
&  \displaystyle \lambda_r \geq 0, & \! \! \! \! \forall r \in \mathcal{R}.\label{stoch:mp:agg:eq4}
\end{eqnarray}
This problem is similar to the one presented in \cite{dantzigmadansky1961} and has $n + 1$ constraints. Let $\overline{x} \in \mathbb{R}^{n}_+$ and $\overline{v} \in \mathbb{R}$ be the dual variables associated with constraints \eqref{stoch:mp:agg:eq1} and \eqref{stoch:mp:agg:eq2}, respectively. We obtain $|\mathcal{S}|$ subproblems, where the $i$-th subproblem takes the form
\begin{equation}
SP^i(\overline{x}) := \max_{\theta_i \in \Theta^i} \left\{(h_i - T_i\overline{x})^T\theta_i\right\}. \label{stoch:subproblem}
\end{equation}
This problem is the dual of \eqref{stoch:prog:sec:obj}-\eqref{stoch:prog:sec:eq2}.
When solving \eqref{stoch:subproblem} we can obtain: (a) an optimal bounded solution or (b) an unbounded solution. 
If all the subproblems lead to an optimal bounded solution, then we obtain an extreme point $\theta_p^i \in \mathbb{R}^{\tilde{m}}$ from each subproblem $SP^i(\overline{x})$. In this case, the reduced cost of the aggregated column is given by 
\begin{equation}
z_{SP}(\overline{x},\overline{v}):= {\max} \{0, \sum_{i \in \mathcal{S}} SP^i(\overline{x}) - \overline{v}\}. \nonumber
\end{equation}
Otherwise, let $\mathcal{U} \subseteq \mathcal{S}$ be the subset of subproblems with unbounded solutions. For each $i \in \mathcal{U}$ we obtain an extreme ray $ \theta_r^i \in \mathbb{R}^{\tilde{m}}$, so the aggregated column is generated by using these extreme rays only. As a result, the reduced cost of the new column is given by
\begin{equation}
z_{SP}(\overline{x},\overline{v}) := {\max} \{0, \sum_{i \in \mathcal{U}} (h_i - T_i\overline{x})^T\theta_r^i \}. \nonumber
\end{equation}

\subsection{Computational experiments}

In this section, we report the results of computational experiments in which we verify the performance of the PDCGM for solving the two-stage stochastic programming problem. We have selected the instances proposed in \cite{ariyawansa2004} and \cite{POSTS}, which have been widely used in the stochastic programming literature. 
All these instances are publicly available in the SMPS format \cite{birge1987}. 
Table \ref{mpc:stochastic:tab:instances} gives the basic information regarding each instance, namely the number of scenarios, the optimal value, and the numbers of columns and rows in the first-stage problem, in the second-stage problem and in the deterministic equivalent problem (DEP). Notice that the same instance name may lead to more than one instance by varying the number of scenarios. From the last two columns of the table, we see that instances with a large number of scenarios have very large-scale corresponding DEPs and hence some of them 
challenge the state-of-the-art implementations of simplex-type methods or interior point methods \cite{zverovichetal2012}. 
On the other hand, these instances can be effectively solved by using a decomposition technique.

\begin{table}[htbp]
\setlength{\tabcolsep}{4.5pt}
\caption{TSSP: statistics for instances from \cite{ariyawansa2004} and \cite{POSTS}}
\scriptsize
\centering
\begin{tabular}{lrrrrrrrr}\hline
 &    &     & \multicolumn{ 2}{c}{Stage 1} & \multicolumn{ 2}{c}{Stage 2} & \multicolumn{ 2}{c}{DEP} \\ \cline{4-5} \cline{6-7} \cline{8-9}
 & $|\mathcal{S}|$ & Optimal value & Cols & Rows & Cols & Rows & Cols & Rows \\ \hline
fxm & 6 & 1.84171E+04 & 92 & 114 & 238 & 343 & 1520 & 2172 \\ 
 & 16 & 1.84168E+04 & 92 & 114 & 238 & 343 & 3900 & 5602 \\ \hline
fxmev & 1 & 1.84168E+04 & 92 & 114 & 238 & 343 & 330 & 457 \\ \hline
pltexpa & 6 & -9.47935E+00 & 62 & 188 & 104 & 272 & 686 & 1820 \\ 
 & 16 & -9.66331E+00 & 62 & 188 & 104 & 272 & 1726 & 4540 \\ \hline
stormg2 & 8 & 1.55352E+07 & 185 & 121 & 528 & 1259 & 4409 & 10193 \\ 
 & 27 & 1.55090E+07 & 185 & 121 & 528 & 1259 & 14441 & 34114 \\ 
 & 125 & 1.55121E+07 & 185 & 121 & 528 & 1259 & 66185 & 157496 \\ 
 & 1000 & 1.58026E+07 & 185 & 121 & 528 & 1259 & 528185 & 1259121 \\ \hline
airl-first & 25 & 2.49102E+05 & 2 & 4 & 6 & 8 & 152 & 204 \\ 
airl-second & 25 & 2.69665E+05 & 2 & 4 & 6 & 8 & 152 & 204 \\ 
airl-randgen & 676 & 2.50262E+05 & 2 & 4 & 6 & 8 & 4058 & 5412 \\ \hline 
assets & 100 & -7.23839E+02 & 5 & 13 & 5 & 13 & 505 & 1313 \\ 
 & 37500 & -6.95963E+02 & 5 & 13 & 5 & 13 & 187505 & 487513 \\ \hline
4node & 1 & 4.13388E+02 & 14 & 52 & 74 & 186 & 88 & 238 \\ 
 & 2 & 4.14013E+02 & 14 & 52 & 74 & 186 & 162 & 424 \\ 
 & 4 & 4.16513E+02 & 14 & 52 & 74 & 186 & 310 & 796 \\ 
 & 8 & 4.18513E+02 & 14 & 52 & 74 & 186 & 606 & 1540 \\ 
 & 16 & 4.23013E+02 & 14 & 52 & 74 & 186 & 1198 & 3028 \\ 
 & 32 & 4.23013E+02 & 14 & 52 & 74 & 186 & 2382 & 6004 \\ 
 & 64 & 4.23013E+02 & 14 & 52 & 74 & 186 & 4750 & 11956 \\ 
 & 128 & 4.23013E+02 & 14 & 52 & 74 & 186 & 9486 & 23860 \\ 
 & 256 & 4.25375E+02 & 14 & 52 & 74 & 186 & 18958 & 47668 \\ 
 & 512 & 4.29963E+02 & 14 & 52 & 74 & 186 & 37902 & 95284 \\ 
 & 1024 & 4.34112E+02 & 14 & 52 & 74 & 186 & 75790 & 190516 \\ 
 & 2048 & 4.41738E+02 & 14 & 52 & 74 & 186 & 151566 & 380980 \\ 
 & 4096 & 4.46856E+02 & 14 & 52 & 74 & 186 & 303118 & 761908 \\ 
 & 8192 & 4.46856E+02 & 14 & 52 & 74 & 186 & 606222 & 1523764 \\ 
 & 16384 & 4.46856E+02 & 14 & 52 & 74 & 186 & 1212430 & 3047476 \\ 
 & 32768 & 4.46856E+02 & 14 & 52 & 74 & 186 & 2424846 & 6094900 \\ \hline
4node-base & 1 & 4.13388E+02 & 16 & 52 & 74 & 186 & 90 & 238 \\ 
 & 2 & 4.14013E+02 & 16 & 52 & 74 & 186 & 164 & 424 \\ 
 & 4 & 4.14388E+02 & 16 & 52 & 74 & 186 & 312 & 796 \\ 
 & 8 & 4.14688E+02 & 16 & 52 & 74 & 186 & 608 & 1540 \\ 
 & 16 & 4.14688E+02 & 16 & 52 & 74 & 186 & 1200 & 3028 \\ 
 & 32 & 4.16600E+02 & 16 & 52 & 74 & 186 & 2384 & 6004 \\ 
 & 64 & 4.16600E+02 & 16 & 52 & 74 & 186 & 4752 & 11956 \\ 
 & 128 & 4.16600E+02 & 16 & 52 & 74 & 186 & 9488 & 23860 \\ 
 & 256 & 4.17162E+02 & 16 & 52 & 74 & 186 & 18960 & 47668 \\ 
 & 512 & 4.20293E+02 & 16 & 52 & 74 & 186 & 37904 & 95284 \\ 
 & 1024 & 4.23050E+02 & 16 & 52 & 74 & 186 & 75792 & 190516 \\ 
 & 2048 & 4.23763E+02 & 16 & 52 & 74 & 186 & 151568 & 380980 \\ 
 & 4096 & 4.24753E+02 & 16 & 52 & 74 & 186 & 303120 & 761908 \\ 
 & 8192 & 4.24775E+02 & 16 & 52 & 74 & 186 & 606224 & 1523764 \\ 
 & 16384 & 4.24775E+02 & 16 & 52 & 74 & 186 & 1212432 & 3047476 \\ 
 & 32768 & 4.24775E+02 & 16 & 52 & 74 & 186 & 2424848 & 6094900 \\ \hline
4node-old & 32 & 8.30941E+04 & 14 & 52 & 74 & 186 & 2382 & 6004 \\ \hline
chem & 2 & -1.30092E+04 & 38 & 39 & 46 & 41 & 130 & 121 \\ 
chem-base & 2 & -1.30092E+04 & 38 & 39 & 40 & 41 & 118 & 121 \\ \hline
lands & 3 & 3.81853E+02 & 2 & 4 & 7 & 12 & 23 & 40 \\ 
lands-blocks & 3 & 3.81853E+02 & 2 & 4 & 7 & 12 & 23 & 40 \\ \hline
env-aggr & 5 & 2.04787E+04 & 48 & 49 & 48 & 49 & 288 & 294 \\ 
env-first & 5 & 1.97774E+04 & 48 & 49 & 48 & 49 & 288 & 294 \\ 
env-loose & 5 & 1.97774E+04 & 48 & 49 & 48 & 49 & 288 & 294 \\ \hline
env & 15 & 2.22653E+04 & 48 & 49 & 48 & 49 & 768 & 784 \\ 
 & 1200 & 2.24289E+04 & 48 & 49 & 48 & 49 & 57648 & 58849 \\ 
 & 1875 & 2.24471E+04 & 48 & 49 & 48 & 49 & 90048 & 91924 \\ 
 & 3780 & 2.24410E+04 & 48 & 49 & 48 & 49 & 181488 & 185269 \\ 
 & 5292 & 2.24384E+04 & 48 & 49 & 48 & 49 & 254064 & 259357 \\ 
 & 8232 & 2.24391E+04 & 48 & 49 & 48 & 49 & 395184 & 403417 \\ 
 & 32928 & 2.24391E+04 & 48 & 49 & 48 & 49 & 1580592 & 1613521 \\ \hline
env-diss-aggr & 5 & 1.59639E+04 & 48 & 49 & 48 & 49 & 288 & 294 \\ 
env-diss-first & 5 & 1.47946E+04 & 48 & 49 & 48 & 49 & 288 & 294 \\ 
env-diss-loose & 5 & 1.47946E+04 & 48 & 49 & 48 & 49 & 288 & 294 \\ \hline
env-diss & 15 & 2.07739E+04 & 48 & 49 & 48 & 49 & 768 & 784 \\ 
 & 1200 & 2.08086E+04 & 48 & 49 & 48 & 49 & 57648 & 58849 \\ 
 & 1875 & 2.08093E+04 & 48 & 49 & 48 & 49 & 90048 & 91924 \\ 
 & 3780 & 2.07947E+04 & 48 & 49 & 48 & 49 & 181488 & 185269 \\ 
 & 5292 & 2.07886E+04 & 48 & 49 & 48 & 49 & 254064 & 259357 \\ 
 & 8232 & 2.07994E+04 & 48 & 49 & 48 & 49 & 395184 & 403417 \\ 
 & 32928 & 2.07994E+04 & 48 & 49 & 48 & 49 & 1580592 & 1613521 \\ \hline
phone1 & 1 & 3.69000E+01 & 1 & 8 & 23 & 85 & 24 & 93 \\ 
phone & 32768 & 3.69000E+01 & 1 & 8 & 23 & 85 & 753665 & 2785288 \\ \hline
stocfor1 & 1 & -4.11320E+04 & 15 & 15 & 102 & 96 & 117 & 111 \\ 
stocfor2 & 64 & -3.97724E+04 & 15 & 15 & 102 & 96 & 6543 & 6159 \\ \hline
\end{tabular}
\label{mpc:stochastic:tab:instances}
\end{table}

In the computational experiments we apply the PDCGM to solve the aggregated master problem formulation (\ref{stoch:mp:agg:obj})-(\ref{stoch:mp:agg:eq4}), using one subproblem of type (\ref{stoch:subproblem}) for each scenario in the instance.  The initial RMP is given by the columns we generate from the first-stage components of the optimal solution of the expected value problem \cite{Birge1997}. In addition, to guarantee the feasibility of any RMP, we have added to it an artificial column with coefficient cost equal to $10^{6}$, in which the entry on the convexity constraint (\ref{stoch:mp:agg:eq2}) is equal to 1 and the remaining entries are equal to 0. The optimality tolerance used to stop the column generation procedure was set to $\delta = 10^{-5}$ and the degree of optimality was set to $D = 5.0$. The experiments were run on a Linux PC with an Intel Core i7 2.8 GHz CPU and 8.0 GB of memory. To solve the subproblems we have used the solver IBM ILOG CPLEX v. 12.4 with its default settings. In Table \ref{mpc:stochastic:tab:results} we present the results of the computational experiments. The first two columns give the problem name and the number of scenarios which are considered in the instance, respectively. Columns 3 and 4 show the number of outer iterations and the total CPU time (in seconds) to solve the instance by using the PDCGM. {Recall that an outer iteration involves solving the RMP and calling the oracle with the aim of generating a new column.} The last row in the table gives the total number of outer iterations and the total CPU time to solve all the instances. Notice that the PDCGM was able to solve all the instances in the set in less than 430 seconds, with the total number of outer iterations equal to 1762. The average CPU time was 5.7 seconds, and the largest CPU time was 81.05 seconds. The number of outer iterations to solve an instance was never larger than 70, and on average it was 23.5. 

\begin{table}[htbp]

\setlength{\tabcolsep}{3pt}
\caption{TSSP: PDCGM results and comparison with the results reported in \cite{zverovichetal2012} for the standard column generation (Benders), and bundle-type method (Level).}
\scriptsize
\centering
\begin{tabular}{lrrrrrrrrrrr}\hline
 & \multicolumn{7}{l}{} & \multicolumn{4}{c}{Relative to PDCGM} \\ \cline{9-12}
 & \multicolumn{1}{l}{} & \multicolumn{ 2}{c}{PDCGM} & \multicolumn{ 2}{c}{Benders} & \multicolumn{ 2}{c}{Level} & \multicolumn{ 2}{c}{Benders} & \multicolumn{ 2}{c}{Level} \\ \cline{3-12}
 & nS & Outer & T (s) & Outer & T (s) & Outer & T (s) & Outer & T (s) & Outer & T (s) \\ \hline
	fxm	&	6	&	16	&	0.21	&	25	&	0.07	&	20	&	0.13	&	1.56	&	0.33	&	1.25	&	0.62	\\	
		&	16	&	18	&	0.25	&	25	&	0.08	&	20	&	0.13	&	1.39	&	0.32	&	1.11	&	0.52	\\	\hline
	fxmev	&	1	&	15	&	0.19	&	25	&	0.07	&	20	&	0.11	&	1.67	&	0.37	&	1.33	&	0.58	\\	\hline
	pltexpa	&	6	&	5	&	0.11	&	1	&	0.02	&	1	&	0.02	&	0.20	&	0.18	&	0.20	&	0.18	\\	
		&	16	&	5	&	0.11	&	1	&	0.02	&	1	&	0.02	&	0.20	&	0.18	&	0.20	&	0.18	\\	\hline
	stormg2	&	8	&	25	&	0.66	&	23	&	0.12	&	20	&	0.14	&	0.92	&	0.18	&	0.80	&	0.21	\\	
		&	27	&	30	&	1.31	&	32	&	0.40	&	17	&	0.26	&	1.07	&	0.31	&	0.57	&	0.20	\\	
		&	125	&	38	&	3.8	&	34	&	1.47	&	17	&	0.79	&	0.89	&	0.39	&	0.45	&	0.21	\\	
		&	1000	&	42	&	22.03	&	41	&	9.83	&	21	&	5.28	&	0.98	&	0.45	&	0.50	&	0.24	\\	\hline
	airl-first	&	25	&	14	&	0.03	&	16	&	0.03	&	17	&	0.03	&	1.14	&	1.00	&	1.21	&	1.00	\\	
	airl-second	&	25	&	12	&	0.05	&	10	&	0.02	&	17	&	0.03	&	0.83	&	0.40	&	1.42	&	0.60	\\	
	airl-randgen	&	676	&	16	&	0.16	&	18	&	0.19	&	18	&	0.19	&	1.13	&	1.19	&	1.13	&	1.19	\\	\hline
	assets	&	100	&	1	&	0.02	&	1	&	0.02	&	1	&	0.03	&	1.00	&	1.00	&	1.00	&	1.50	\\	
		&	37500	&	1	&	1.28	&	2	&	74.53	&	2	&	74.42	&	2.00	&	58.23	&	2.00	&	58.14	\\	\hline
	4node	&	1	&	29	&	0.18	&	24	&	0.03	&	21	&	0.05	&	0.83	&	0.17	&	0.72	&	0.28	\\	
		&	2	&	31	&	0.23	&	38	&	0.03	&	42	&	0.09	&	1.23	&	0.13	&	1.35	&	0.39	\\	
		&	4	&	31	&	0.23	&	41	&	0.03	&	45	&	0.09	&	1.32	&	0.13	&	1.45	&	0.39	\\	
		&	8	&	32	&	0.26	&	64	&	0.06	&	45	&	0.09	&	2.00	&	0.23	&	1.41	&	0.35	\\	
		&	16	&	37	&	0.32	&	67	&	0.09	&	44	&	0.13	&	1.81	&	0.28	&	1.19	&	0.41	\\	
		&	32	&	37	&	0.38	&	100	&	0.20	&	51	&	0.19	&	2.70	&	0.53	&	1.38	&	0.50	\\	
		&	64	&	36	&	0.45	&	80	&	0.23	&	54	&	0.31	&	2.22	&	0.51	&	1.50	&	0.69	\\	
		&	128	&	37	&	0.56	&	74	&	0.33	&	50	&	0.40	&	2.00	&	0.59	&	1.35	&	0.71	\\	
		&	256	&	35	&	0.86	&	71	&	0.81	&	48	&	0.74	&	2.03	&	0.94	&	1.37	&	0.86	\\	
		&	512	&	40	&	2.06	&	92	&	3.16	&	51	&	1.80	&	2.30	&	1.53	&	1.28	&	0.87	\\	
		&	1024	&	38	&	3.56	&	70	&	4.37	&	53	&	3.36	&	1.84	&	1.23	&	1.39	&	0.94	\\	
		&	2048	&	37	&	6.02	&	83	&	10.01	&	49	&	6.65	&	2.24	&	1.66	&	1.32	&	1.10	\\	
		&	4096	&	43	&	10.83	&	89	&	15.69	&	46	&	7.75	&	2.07	&	1.45	&	1.07	&	0.72	\\	
		&	8192	&	48	&	22.55	&	106	&	39.58	&	55	&	19.28	&	2.21	&	1.76	&	1.15	&	0.85	\\	
		&	16384	&	42	&	37.97	&	110	&	84.15	&	52	&	38.45	&	2.62	&	2.22	&	1.24	&	1.01	\\	
		&	32768	&	40	&	54.14	&	122	&	165.48	&	62	&	108.68	&	3.05	&	3.06	&	1.55	&	2.01	\\	\hline
	4node-base	&	1	&	22	&	0.13	&	31	&	0.03	&	16	&	0.03	&	1.41	&	0.23	&	0.73	&	0.23	\\	
		&	2	&	28	&	0.2	&	44	&	0.03	&	29	&	0.05	&	1.57	&	0.15	&	1.04	&	0.25	\\	
		&	4	&	29	&	0.2	&	58	&	0.05	&	30	&	0.06	&	2.00	&	0.25	&	1.03	&	0.30	\\	
		&	8	&	28	&	0.2	&	47	&	0.04	&	35	&	0.09	&	1.68	&	0.20	&	1.25	&	0.45	\\	
		&	16	&	30	&	0.27	&	56	&	0.07	&	30	&	0.09	&	1.87	&	0.26	&	1.00	&	0.33	\\	
		&	32	&	35	&	0.37	&	63	&	0.14	&	37	&	0.14	&	1.80	&	0.38	&	1.06	&	0.38	\\	
		&	64	&	35	&	0.44	&	61	&	0.20	&	33	&	0.19	&	1.74	&	0.45	&	0.94	&	0.43	\\	
		&	128	&	33	&	0.47	&	65	&	0.33	&	37	&	0.30	&	1.97	&	0.70	&	1.12	&	0.64	\\	
		&	256	&	39	&	0.9	&	66	&	0.76	&	31	&	0.45	&	1.69	&	0.84	&	0.79	&	0.50	\\	
		&	512	&	38	&	1.7	&	84	&	2.78	&	37	&	1.23	&	2.21	&	1.64	&	0.97	&	0.72	\\	
		&	1024	&	49	&	4.37	&	115	&	8.13	&	41	&	2.83	&	2.35	&	1.86	&	0.84	&	0.65	\\	
		&	2048	&	55	&	7.79	&	142	&	16.76	&	42	&	5.21	&	2.58	&	2.15	&	0.76	&	0.67	\\	
		&	4096	&	59	&	13.46	&	174	&	32.73	&	39	&	9.01	&	2.95	&	2.43	&	0.66	&	0.67	\\	
		&	8192	&	70	&	31.11	&	290	&	113.43	&	48	&	21.24	&	4.14	&	3.65	&	0.69	&	0.68	\\	
		&	16384	&	68	&	59.55	&	175	&	139.46	&	41	&	40.21	&	2.57	&	2.34	&	0.60	&	0.68	\\	
		&	32768	&	61	&	81.05	&	191	&	267.16	&	49	&	86.95	&	3.13	&	3.30	&	0.80	&	1.07	\\	\hline
	4node-old	&	32	&	21	&	0.4	&	30	&	0.07	&	20	&	0.08	&	1.43	&	0.18	&	0.95	&	0.20	\\	\hline
	chem	&	2	&	5	&	0.02	&	7	&	0.03	&	15	&	0.03	&	1.40	&	1.50	&	3.00	&	1.50	\\	
	chem-base	&	2	&	10	&	0.04	&	6	&	0.02	&	14	&	0.04	&	0.60	&	0.50	&	1.40	&	1.00	\\	\hline
	lands	&	3	&	7	&	0.03	&	8	&	0.02	&	10	&	0.02	&	1.14	&	0.67	&	1.43	&	0.67	\\	
	lands-blocks	&	3	&	7	&	0.01	&	8	&	0.01	&	10	&	0.02	&	1.14	&	1.00	&	1.43	&	2.00	\\	\hline
	env-aggr	&	5	&	5	&	0.05	&	3	&	0.02	&	16	&	0.03	&	0.60	&	0.40	&	3.20	&	0.60	\\	
	env-first	&	5	&	2	&	0.03	&	1	&	0.02	&	1	&	0.02	&	0.50	&	0.67	&	0.50	&	0.67	\\	
	env-loose	&	5	&	2	&	0.03	&	1	&	0.01	&	1	&	0.02	&	0.50	&	0.33	&	0.50	&	0.67	\\	\hline
	env	&	15	&	5	&	0.05	&	3	&	0.03	&	15	&	0.04	&	0.60	&	0.60	&	3.00	&	0.80	\\	
		&	1200	&	5	&	0.45	&	3	&	0.29	&	15	&	1.47	&	0.60	&	0.64	&	3.00	&	3.27	\\	
		&	1875	&	5	&	0.68	&	3	&	0.48	&	15	&	2.38	&	0.60	&	0.71	&	3.00	&	3.50	\\	
		&	3780	&	5	&	1.2	&	3	&	1.07	&	15	&	4.65	&	0.60	&	0.89	&	3.00	&	3.88	\\	
		&	5292	&	5	&	1.6	&	3	&	1.67	&	15	&	6.52	&	0.60	&	1.04	&	3.00	&	4.08	\\	
		&	8232	&	5	&	2.41	&	3	&	3.15	&	15	&	10.69	&	0.60	&	1.31	&	3.00	&	4.44	\\	
		&	32928	&	5	&	9.64	&	3	&	33.90	&	15	&	64.32	&	0.60	&	3.52	&	3.00	&	6.67	\\	\hline
	env-diss-aggr	&	5	&	11	&	0.17	&	9	&	0.03	&	22	&	0.04	&	0.82	&	0.18	&	2.00	&	0.24	\\	
	env-diss-first	&	5	&	9	&	0.06	&	14	&	0.02	&	12	&	0.03	&	1.56	&	0.33	&	1.33	&	0.50	\\	
	env-diss-loose	&	5	&	9	&	0.05	&	15	&	0.03	&	5	&	0.03	&	1.67	&	0.60	&	0.56	&	0.60	\\	\hline
	env-diss	&	15	&	17	&	0.28	&	27	&	0.04	&	35	&	0.09	&	1.59	&	0.14	&	2.06	&	0.32	\\	
		&	1200	&	10	&	0.82	&	24	&	0.96	&	35	&	2.38	&	2.40	&	1.17	&	3.50	&	2.90	\\	
		&	1875	&	18	&	1.93	&	29	&	2.13	&	36	&	3.82	&	1.61	&	1.10	&	2.00	&	1.98	\\	
		&	3780	&	13	&	2.49	&	29	&	4.28	&	36	&	7.54	&	2.23	&	1.72	&	2.77	&	3.03	\\	
		&	5292	&	17	&	4.34	&	34	&	6.92	&	38	&	11.01	&	2.00	&	1.59	&	2.24	&	2.54	\\	
		&	8232	&	16	&	5.86	&	35	&	12.08	&	41	&	19.12	&	2.19	&	2.06	&	2.56	&	3.26	\\	
		&	32928	&	16	&	22.37	&	35	&	67.59	&	41	&	95.59	&	2.19	&	3.02	&	2.56	&	4.27	\\	\hline
	phone1	&	1	&	1	&	0.01	&	1	&	0.02	&	1	&	0.02	&	1.00	&	2.00	&	1.00	&	2.00	\\	
	phone	&	32768	&	1	&	0.76	&	1	&	41.09	&	1	&	41.00	&	1.00	&	54.07	&	1.00	&	53.95	\\	\hline
	stocfor1	&	1	&	11	&	0.11	&	6	&	0.02	&	6	&	0.03	&	0.55	&	0.18	&	0.55	&	0.27	\\	
	stocfor2	&	64	&	9	&	0.11	&	7	&	0.09	&	9	&	0.10	&	0.78	&	0.82	&	1.00	&	0.91	\\	\hline
\bf	Total	&		&	1762	&	429.02	&	3498	&	1169.31	&	2045	&	708.85	&	1.99	&	2.73	&	1.16	&	1.65	\\	
\hline
\end{tabular}
\label{mpc:stochastic:tab:results}
\end{table}

To verify the performance of the PDCGM in relation to other column generation/cutting plane approaches, we have added to Table \ref{mpc:stochastic:tab:results} the results presented in \cite{zverovichetal2012} for the instances given in Table \ref{mpc:stochastic:tab:instances}. More specifically, we borrow from \cite{zverovichetal2012} the results which are reported for the standard Benders decomposition \cite{Ben62,vanslyke1969} and the Level-set method \cite{lemarechal1995}. From a column generation point of view, these two methods correspond to solving the aggregated master problem formulation (\ref{stoch:mp:agg:obj})-(\ref{stoch:mp:agg:eq4}) by the standard column generation and by a bundle-type method (a stabilized column generation variant), respectively. The number of outer iterations and CPU time (in seconds) for solving the instances by these two approaches are given in columns 5 to 8 of Table \ref{mpc:stochastic:tab:results}. The remaining columns in the table show the relative performance of these methods in relation to PDCGM, \textit{i.e.}, the ratio between the values in columns 5 to 8 and the corresponding values in columns 3 and 4. {We have taken the iteration numbers exactly as they are presented in \cite{zverovichetal2012}, while the CPU times have been scaled according to the benchmark available at \url{https://www.cpubenchmark.net/singleThread.html}. The computer used in \cite{zverovichetal2012} was a Linux PC with an Intel Core i5 2.4 GHz CPU and 6.0 GB of memory, and has a score of 1355 whereas the machine we used to run the PDCGM experiments has a score of 1602. Hence, we have multiplied their CPU times by a factor of $0.85$. Moreover, the authors have implemented the methods on top of the FortSP stochastic solver system \cite{ellison2010}, a state-of-the-art solver for stochastic programming. Therefore, the conclusions about CPU times should be taken \emph{cautiously} and, hence, we focus on the number of outer iterations. }

The results in Table \ref{mpc:stochastic:tab:results} show that even though the standard column generation (Benders) had the smallest CPU times on instances with few scenarios, this method delivers the worst overall performance in relation to the PDCGM and the bundle-type method (Level). To solve all the instances, the standard approach required 3498 iterations, { which is almost twice the figures obtained by the PDCGM. 
This is justified by the good performance of the PDCGM on the instances with a large number of scenarios, an important feature in practice. Besides, the overall performance of the PDCGM was similar to that of the level method, although the latter had about 16\% more iterations in total. Hence, the results indicate that the PDCGM is competitive with respect to the level method, which is considered an efficient method for solving two-stage stochastic programming problems \cite{zverovichetal2012}. }

\section{Multicommodity network flow problem (MCNF)} \label{sec:nlmcnf}

Multicommodity network flow problems (MCNF) have been widely studied in the literature and can be applied in contexts in which commodities (e.g., goods and data packages) must be transported through a network with limited capacity and arc costs \cite{ouorou2000}. The current real-life applications involve transportation as well as telecommunication networks with a large number of arcs and commodities. Hence they lead to very large-scale optimization problems which require efficient solution methods. 
Column generation approaches have been successfully used for solving this class of problems \cite{GofGonSarVia96,babonneau2006,lemarechal2009,babonneau2009,frangioni2013stabilized}. 
{ In this work, we consider one of the most basic variants of the MCNF problem which includes a linear cost function and where the cost depends on the flow traversing an arc and it is independent of the commodity.}
In this section, we first describe the column generation formulation of the MCNF with linear objective costs. Then, we present the results of using the primal-dual column generation technique for solving large-scale instances, some of them taken from real-life applications. 

\subsection{Problem formulation}

Consider a set $\mathcal{K} = \{ 1, \ldots, K\}$ of commodities which must be routed through a given network represented by the set of nodes $\mathcal{N} = \{ 1, \ldots, n\}$ and the set of arcs $\mathcal{M} = \{ 1, \ldots, m\}$. For each commodity $k \in \mathcal{K}$ there is a source node $s_k \in \mathcal{N}$ and a sink node $t_k \in \mathcal{N}$, so that the total demand of the commodity ($d_k$) must be routed from $s_k$ to $t_k$ using one or more arcs of the underlying network. Let $A$ be the $n \times m$ node-arc incidence matrix of the network determined by sets $\mathcal{N}$ and $\mathcal{M}$. In order to associate the demand of each commodity $k$ with all the nodes in the network, we define an $n$-vector $b^{k}$ as follows
\begin{equation}
b^{k}_{i} = \left\{ 
   \begin{array}{rl}
   d_k, & \mbox{ if } i = s_k, \\ 
   -d_k, & \mbox{ if } i = t_k, \\
   0, & \mbox{ otherwise}, 
   \end{array} \nonumber
\right. \ \ \ \ \ \forall i \in \mathcal{N}. \nonumber
\end{equation}
Let $x^{k}_{ij}$ be the decision variable that determines the flow of commodity $k \in \mathcal{K}$ assigned to arc $(i,j) \in \mathcal{M}$. The total flow assigned to a given arc $(i,j) \in \mathcal{M}$ 
cannot exceed
the arc capacity $C_{ij}$. In addition, there is a cost $t_{ij}$ that depends linearly on the total flow assigned to the arc. A compact formulation of the (linear) MCNF is given by
\begin{eqnarray}
\min_{x} & \displaystyle {\sum_{k \in \mathcal{K}}} \sum_{(i,j) \in \mathcal{M}} t_{ij} x^{k}_{ij}, & \label{eq:lmcnf:obj}\\
\mbox{s.t.} & \displaystyle \sum_{k \in \mathcal{K}} x^{k}_{ij} \leq C_{ij}, & \ \ \forall (i,j) \in \mathcal{M}, \label{eq:lmcnf:01}\\
& \displaystyle A x^{k} = b^{k}, & \ \ \forall k \in \mathcal{K}, \label{eq:lmcnf:02}\\
& \displaystyle x^{k}_{ij} \geq 0, & \ \ \forall k \in \mathcal{K}, \ \forall (i,j) \in \mathcal{M}. \label{eq:lmcnf:03}
\end{eqnarray}
This formulation typically leads to a large-scale problem when modeling real-life situations with many commodities and arcs in the network, as the number of variables and constraints may become very large. Solving \eqref{eq:lmcnf:obj}-\eqref{eq:lmcnf:03} by a linear optimization method may be prohibitive in practice, even for the current state-of-the-art solvers.
Fortunately, the coefficient matrix of this formulation has a special structure that can be exploited to obtain a more efficient solution strategy.

In this paper, we focus on MCNF problems in which the costs in the network depend on the arcs only and the commodities compete equally for the capacity on the arcs. These problems are usually solved by approaches which are based on column generation procedures \cite{lemarechal2009,babonneau2006,babonneau2009}. As recognized in \cite{babonneau2006}, there is a different type of MCNF problems in which the costs depend additionally on the commodities assigned to the arc, and the commodities compete for mutual and/or individual capacities \cite{frangioni1999,castro2003,castro2012}. No paper in the MCNF literature deals with both types of problems simultaneously.

\subsection{Decomposition and column generation formulation}\label{mpcpaper:mcnf:lmcnf}

In this section, we apply the Dantzig-Wolfe decomposition (DWD) to the compact formulation of the MCNF, following the description given in Section \ref{sec:dwdcolgen}. The coefficient matrix of formulation (\ref{eq:lmcnf:obj})-(\ref{eq:lmcnf:03}) has a special structure. Namely, the incidence matrix $A$ is replicated $K$ times resulting in a block-angular structure, making appropriate the use of the DWD. Since the decision variables associated with different commodities are only connected by the linking constraints (\ref{eq:lmcnf:01}), we define the following $K$ independent subsets
\begin{equation}
\mathcal{X}^{k} = \{ x^{k} \mid  A x^{k} = b^{k}, \ x^{k}_{ij} \geq 0, \ \forall (i,j) \in \mathcal{M} \}, \ \ k \in \mathcal{K}. \label{eq:lmcnf:indepsubsets}
\end{equation}
Let ${P}_k$ be the set of indices of extreme points in $\mathcal{X}^{k}$, $k \in \mathcal{K}$. Each set $\mathcal{X}^{k}$ is bounded and hence it can be fully represented by its extreme points $x^{k}_{p}$, $p \in P_k$. As a result, any point $x^{k} \in \mathcal{X}^{k}$ can be rewritten similarly to (\ref{eq:rewritex}). 
The master problem related to the compact formulation (\ref{eq:lmcnf:obj})-(\ref{eq:lmcnf:03}) is then given by
\begin{eqnarray}
\min_{\lambda} & \displaystyle {\sum_{k \in \mathcal{K}}} \sum_{p \in P_k} \sum_{(i,j) \in \mathcal{M}} t_{ij} (x^{k}_{p})_{ij} \lambda^{k}_p, & \label{eq:lmcnf:mp:obj}\\
\mbox{s.t.} & \displaystyle {\sum_{k \in \mathcal{K}}} \sum_{p \in P_k} (x^{k}_{p})_{ij} \lambda^{k}_p \leq C_{ij}, & \ \ \forall (i,j) \in \mathcal{M}, \label{eq:lmcnf:mp:01}\\
& \displaystyle \sum_{p \in P_k} \lambda^{k}_p = 1, & \ \ \forall k \in \mathcal{K},\label{eq:lmcnf:mp:02} \\
& \displaystyle \lambda^{k}_{p} \geq 0, & \ \ \forall k \in \mathcal{K}, \ {\forall}p \in P_k. \label{eq:lmcnf:mp:03}
\end{eqnarray}
Since the number of variables $\lambda^{k}_{p}$ can be huge, we recur to the column generation technique for solving the master problem, as described in Section \ref{sec:dwdcolgen}.
The corresponding oracle is given by $K$ subproblems, so that the $k$-th subproblem provides an extreme point of $\mathcal{X}^{k}$ which consists in the shortest {reduced cost} path that routes the total demand of commodity $k$ from source to sink nodes. 
Let $\overline{u} \in \mathbb{R}^m_{-}$ and $\overline{v} \in \mathbb{R}^K$ be the dual solution vectors
associated with constraints (\ref{eq:lmcnf:mp:01}) and (\ref{eq:lmcnf:mp:02}), respectively. 
The $k$-th subproblem can be stated as
\begin{equation}
{SP}^{k}(\overline{u}) := \min_{x^{k} \in \mathcal{X}^k} \left\{\sum_{(i,j) \in \mathcal{M}}(t_{ij} - \overline{u}_{ij}) x^{k}_{ij} \right\}. \label{eq:lmcnf:spk}
\end{equation}
It is a shortest path problem with arc lengths $t_{ij} - \overline{u}_{ij}$, initial node $s_k$ and end node $t_k$. Since $t_{ij} - \overline{u}_{ij} \geq 0$ for all $(i,j) \in \mathcal{M}$, these subproblems can be solved by the Dijkstra algorithm \cite{dijkstra1959}. 
From the optimal solution of each subproblem, we generate a new column with {reduced} cost ${z}^{k}_{SP}(\overline{u}, \overline{v}):= \min \{0;SP^k(\overline{u})-\overline{v}_k\}$
and the value of the oracle is
\begin{equation}
{z}_{SP}(\overline{u}, \overline{v}):= \sum_{k \in \mathcal{K}} z^k_{SP}(\overline{u}, \overline{v}).
\end{equation}

\subsection{Computational experiments}

In this section, we verify the performance of the PDCGM to solve the master problem formulation of the MCNF. We have run computational experiments based on two sets of instances that have been proposed in \cite{larsson2004} and, since then, they have been used as a benchmark for linear and nonlinear MCNF problems \cite{babonneau2006,alvelos2007,babonneau2009}. The instances in the first set simulate telecommunication networks which are represented by planar graphs. The instances in the second set represent networks with a grid structure and, hence, the number of different paths between two nodes can be very large. All instances are publicly available and can be downloaded from \url{http://www.di.unipi.it/di/groups/optimize/Data/MMCF.html}.
The initial columns in the RMP were generated by using the solution we obtain from calling the oracle with $\overline{u} = 0$. The tolerance of the column generation was set to $\delta = 10^{-5}$ and the degree of
optimality was set to $D = 10$. The experiments were run on a Linux PC with an Intel Core i7 2.8 GHz CPU and 8.0 GB of memory.

In practical MCNF problems, it is very likely that in the optimal solution only a small number of arcs will have full capacity usage. 
As a consequence, many of the arc capacity constraints (\ref{eq:lmcnf:mp:01}) will be inactive at the optimum. If it was possible to identify these constraints before 
solving the problem,
they could be discarded 
without changing the optimal solution. Since knowing the inactive constraints in advance is not possible in general, an active set strategy is used during the optimization process, in order to guess which arc capacity constraints should be classified as active. This idea has been successfully applied in the context of MCNF problems \cite{mcbride1998,frangioni1999,babonneau2006} and, hence, we have added this in our implementation as well. 
In the first RMP, we assume all the arc capacity constraints are inactive, i.e., they are excluded from the formulation.
After each call to the oracle, we verify for each constraint, if the total flow on the corresponding arc violates its capacity. 
If so, the constraint is added to the active set and included in the formulation.
An active constraint can also become inactive if the total flow on the corresponding arc is smaller than a fraction $\gamma$ of the capacity of the arc. In the computational experiments presented below, we have adopted $\gamma = 0.9$.

In Table \ref{mpc:mcnf:tab:resuls1}, columns 1 to 4 show the name of the instances, the number of nodes ($n$), the number of arcs ($m$), and the number of commodities ($K$), respectively. Recall that the number of constraints in the RMP is $m$ + $K$ (disaggregated formulation). The results obtained by PDCGM are presented in the remaining columns. Column 5 shows the optimal value obtained with the optimality tolerance $\delta = 10^{-5}$. Columns 6 to 8 show, the percentage of arc capacity constraints which are active at optimum (\% Act),
the number of outer iterations (Outer) and the total CPU time in seconds for solving the problem (CPU), respectively. The last column gives the percentage of the total CPU time which was used to solve the oracle subproblems (\%Or).

According to the results in Table \ref{mpc:mcnf:tab:resuls1}, the PDCGM was able to solve every instance in the planar set in less than 7103 seconds, and every instance in the grid set in less than 73 seconds.
Most of the CPU time is spent on solving the RMPs, a typical behavior when disaggregated formulations are used for the MCNF \cite{GofGonSarVia96}. The number of outer iterations was on average 22.4 for the planar instances and 20.1 for the grid instances. The number of active arcs was relatively small and never larger than 30\% of the total number of arcs. 

\begin{table}[htbp]
\caption{MCNF: dimensions and computational results obtained by PDCGM.}
\scriptsize
\centering
\begin{tabular}{lcccccccc}
\hline \noalign{\smallskip}
\multicolumn{1}{c}{Instance} & $n$ & $m$ & $K$ & Optimal & \%Act  & Outer & CPU (s) & \%Or \\ 
\noalign{\smallskip} \hline \noalign{\smallskip} 
planar30 & 30 & 150 & 92 & 4.43508E+07 & 12.0 & 13 & 0.04 & 2.3 \\ 
planar50 & 50 & 250 & 267 & 1.22200E+08 & 12.8 & 17 & 0.14 & 0.7 \\ 
planar80 & 80 & 440 & 543 & 1.82438E+08 & 25.7 & 20 & 0.55 & 3.8 \\ 
planar100 & 100 & 532 & 1085 & 2.31340E+08 & 16.5 & 18 & 0.82 & 2.4 \\ 
planar150 & 150 & 850 & 2239 & 5.48087E+08 & 29.4 & 26 & 4.81 & 1.8 \\ 
planar300 & 300 & 1680 & 3584 & 6.89979E+08 & 8.6 & 18 & 4.97 & 5.4 \\ 
planar500 & 500 & 2842 & 3525 & 4.81983E+08 & 2.3 & 13 & 3.32 & 23.7 \\ 
planar800 & 800 & 4388 & 12756 & 1.16737E+08 & 3.4 & 23 & 32.47 & 16.7 \\ 
planar1000 & 1000 & 5200 & 20026 & 3.44962E+09 & 10.8 & 28 & 156.25 & 8.3 \\ 
planar2500 & 2500 & 12990 & 81430 & 1.26623E+10 & 15.0 & 48 & 7102.27 & 4.8 \\ 
\noalign{\smallskip} 
grid1 & 25 & 80 & 50 & 8.27319E+05 & 10.0 & 12 & 0.03 & 4.0 \\ 
grid2 & 25 & 80 & 100 & 1.70538E+06 & 30.0 & 14 & 0.04 & 2.6 \\ 
grid3 & 100 & 360 & 50 & 1.52464E+06 & 4.4 & 13 & 0.05 & 2.0 \\ 
grid4 & 100 & 360 & 100 & 3.03170E+06 & 9.7 & 19 & 0.09 & 4.4 \\ 
grid5 & 225 & 840 & 100 & 5.04969E+06 & 5.4 & 20 & 0.19 & 13.9 \\ 
grid6 & 225 & 840 & 200 & 1.04018E+07 & 15.8 & 24 & 0.49 & 10.7 \\ 
grid7 & 400 & 1520 & 400 & 2.58641E+07 & 8.9 & 21 & 0.92 & 28.9 \\ 
grid8 & 625 & 2400 & 500 & 4.17114E+07 & 13.2 & 28 & 4.28 & 23.8 \\ 
grid9 & 625 & 2400 & 1000 & 8.26531E+07 & 18.7 & 33 & 11.30 & 15.2 \\ 
grid10 & 625 & 2400 & 2000 & 1.64111E+08 & 18.1 & 33 & 13.19 & 15.1 \\ 
grid11 & 625 & 2400 & 4000 & 3.29259E+08 & 12.0 & 15 & 7.77 & 12.9 \\ 
grid12 & 900 & 3480 & 6000 & 5.77187E+08 & 7.0 & 14 & 11.46 & 23.0 \\ 
grid13 & 900 & 3480 & 12000 & 1.15932E+09 & 9.1 & 18 & 30.83 & 11.0 \\ 
grid14 & 1225 & 4760 & 16000 & 1.80268E+09 & 4.0 & 16 & 32.93 & 22.0 \\ 
grid15 & 1225 & 4760 & 32000 & 3.59352E+09 & 4.7 & 22 & 72.90 & 13.6 \\ 
\noalign{\smallskip} \hline
\end{tabular}
\label{mpc:mcnf:tab:resuls1}
\end{table}

To verify the performance of the PDCGM in relation to the best results available in the literature, we have selected one of the most efficient methods for solving the linear MCNF. In \cite{babonneau2006}, the authors use the analytic center cutting plane method (ACCPM) to solve the dual problem of an aggregated master problem formulation of the MCNF. In addition, they use an active set strategy and the elimination of columns. Table \ref{mpc:mcnf:tab:resuls2} shows the best results presented in \cite{babonneau2006} for solving the same instances described in Table \ref{mpc:mcnf:tab:resuls1}. The first column gives the name of the instances.
Columns 2 to 5 show the percentage of arc capacity constraints which are active when the column generation terminates (\%Act), the number of outer iterations (Outer), the total CPU time in seconds, and the percentage of the total CPU time that was spent in the oracle (\%Or), respectively, as reported in \cite{babonneau2006} for the ACCPM-based approach. 
{We have scaled the CPU times using a factor of $0.4$, as these other results were obtained on a Linux PC with Intel Pentium IV 2.8 Ghz and 2 GB of RAM. According to the benchmark provided at \url{https://www.cpubenchmark.net/singleThread.html}, this machine has a score of 635 whereas the machine used to run the PDCGM results has a score of 1602.}
Columns 6 to 9 show the ratios between the results presented in columns 2 to 5 and the corresponding results for PDCGM presented in Table \ref{mpc:mcnf:tab:resuls1}. 
{
Table \ref{mpc:mcnf:tab:resuls2} has the purpose of giving an idea of the overall performance of the PDCGM (with its best settings) in relation to the best results presented in the literature for the linear MCNF. These results are just informative, as different computational environments were used and an aggregated formulation was used in \cite{babonneau2006}. From our experience, we observed that using a disaggregated approach reduces considerably the number of column generation iterations for the MCNF, although the RMP grows quickly and may become more difficult to solve. Nevertheless, a reduction in CPU times can be observed for most of the instances, which indicates that the PDCGM is an efficient alternative for solving the linear MCNF.}

\begin{table}[htbp]

\setlength{\tabcolsep}{6pt}
\caption{MCNF: comparison between PDCGM results and the best results reported in \cite{babonneau2006} for an ACCPM-based approach.}
\scriptsize
\centering
\begin{threeparttable}
\begin{tabular}{lcccccccc}
\hline\noalign{\smallskip}
Instance & \multicolumn{ 4}{c}{ACCPM-based approach\cite{babonneau2006}} & \multicolumn{ 4}{c}{ACCPM/PDCGM} \\ 
 & \multicolumn{1}{l}{\%Act} & \multicolumn{1}{c}{Outer} & \multicolumn{1}{c}{CPU (s)} & \multicolumn{1}{c}{\%Or} & \multicolumn{1}{c}{\%Act} & \multicolumn{1}{c}{Outer} & \multicolumn{1}{c}{CPU (s)} & \multicolumn{1}{c}{\%Or} \\ 
\noalign{\smallskip}\hline\noalign{\smallskip}
planar30	&	12.7	&	48	&	0.16	&	29	&	1.06	&	3.69	&	4.00	&	12.61	\\
planar50	&	13.2	&	97	&	0.44	&	32	&	1.03	&	5.71	&	3.14	&	45.71	\\
planar80	&	25.7	&	283	&	2.6	&	28	&	1.00	&	14.15	&	4.73	&	7.37	\\
planar100	&	17.5	&	257	&	2.08	&	34	&	1.06	&	14.28	&	2.54	&	14.17	\\
planar150	&	29	&	820	&	25.8	&	13	&	0.99	&	31.54	&	5.36	&	7.22	\\
planar300	&	9.1	&	325	&	3.88	&	27	&	1.06	&	18.06	&	0.78	&	5.00	\\
planar500	&	2.6	&	118	&	4.2	&	90	&	1.13	&	9.08	&	1.27	&	3.80	\\
planar800	&	3.6	&	252	&	24.28	&	91	&	1.06	&	10.96	&	0.75	&	5.45	\\
planar1000	&	10.4	&	890	&	229.04	&	55	&	0.96	&	31.79	&	1.47	&	6.63	\\
planar2500	&	15.8	&	3009	&	11782.92	&	28	&	1.05	&	62.69	&	1.66	&	5.83	\\
\noalign{\smallskip}		
grid1	&	12.5	&	24	&	0.08	&	25	&	1.25	&	2.00	&	2.67	&	6.25	\\
grid2	&	32.5	&	52	&	0.24	&	31	&	1.08	&	3.71	&	6.00	&	11.92	\\
grid3	&	5	&	32	&	0.12	&	39	&	1.14	&	2.46	&	2.40	&	19.50	\\
grid4	&	10.3	&	66	&	0.28	&	46	&	1.06	&	3.47	&	3.11	&	10.45	\\
grid5	&	6	&	75	&	0.64	&	58	&	1.11	&	3.75	&	3.37	&	4.17	\\
grid6	&	20.6	&	239	&	3.32	&	46	&	1.30	&	9.96	&	6.78	&	4.30	\\
grid7	&	9.3	&	228	&	5.48	&	70	&	1.04	&	10.86	&	5.96	&	2.42	\\
grid8	&	13.2	&	528	&	39	&	50	&	1.00	&	18.86	&	9.11	&	2.10	\\
grid9	&	18.1	&	720	&	85.12	&	38	&	0.97	&	21.82	&	7.53	&	2.50	\\
grid10	&	18	&	722	&	86.24	&	41	&	0.99	&	21.88	&	6.54	&	2.72	\\
grid11	&	12.2	&	458	&	33.96	&	62	&	1.02	&	30.53	&	4.37	&	4.81	\\
grid12	&	7.4	&	329	&	35.56	&	84	&	1.06	&	23.50	&	3.10	&	3.65	\\
grid13	&	9.6	&	460	&	54.72	&	74	&	1.05	&	25.56	&	1.77	&	6.73	\\
grid14	&	4.4	&	252	&	42.8	&	94	&	1.10	&	15.75	&	1.30	&	4.27	\\
grid15	&	4.9	&	294	&	52.68	&	91	&	1.04	&	13.36	&	0.72	&	6.69	\\
\noalign{\smallskip} \hline
\end{tabular}
\end{threeparttable}
\label{mpc:mcnf:tab:resuls2}
\end{table}

\section{Conclusions}
\label{sec:conclusions}
In this paper we have presented computational evidence of the performance of the primal-dual column generation method (PDCGM) {for} solving the multiple kernel learning (MKL), the two-stage stochastic programming (TSSP) and the multicommodity network flow (MCNF) problems. 
{
We have demonstrated how the Dantzig-Wolfe decomposition and the column generation algorithm can be applied to general convex optimization problems, including the use of aggregated and disaggregated formulations. Additionally, we have provided a thorough presentation of the key components involved in 
applying these techniques to each addressed application.
The performance of the PDCGM was compared against some of the best results available in the literature for MKL, TSSP and MCNF, using different types of column generation/cutting plane methods. }
These applications provided us with different conditions to test the PDCGM, namely disaggregated and aggregated master problem formulations as well as bounded and unbounded subproblems.
{
The computational experiments presented in this paper provide 
extensive evidence that the PDCGM is suitable and competitive in a broader context of optimization than that previously addressed in \cite{GonGonMun2013,MunGon2012}.
}
It is worth mentioning that the PDCGM software is freely available for research use and can be downloaded at \url{http://www.maths.ed.ac.uk/~gondzio/‎software/pdcgm.html}.

Further studies will involve extending the PDCGM to solve a wider range of optimization problems, including those which are defined by nonconvex functions.
In addition, handling second-order oracle information, as in \cite{babonneau2009}, seems to be a promising avenue for future research for the PDCGM.

\section*{Acknowledgements}
We would like to express our gratitude to Victor Zverovich for kindly making available to us some of the TSSP instances included in this study.
Also, we would like to thank Robert Gower for proofreading an early version of this paper.
Also, we are very thankful to the anonymous referees for their careful reading and the important suggestions made, which certainly helped on improving the first draft of this paper.
Pablo Gonz{\'a}lez-Brevis has been partially supported by FONDECYT, Chile through grant 11140521. Pedro Munari has been supported by FAPESP (S\~{a}o Paulo Research Foundation, Brazil) through grant 14/00939-8.

\bibliographystyle{spmpsci}      
\bibliography{mpcbib}   

\end{document}